\documentclass[preprint,3p,12pt]{elsarticle}
\usepackage{amsmath}
\usepackage{stmaryrd}
\usepackage{bbding}
\usepackage{dcolumn}
\usepackage{graphicx}
\usepackage{amsfonts}
\usepackage{amssymb}
\usepackage{psfrag}
\usepackage{wrapfig}
\usepackage{subfigure}
\usepackage{makeidx}
\usepackage{bm}
\usepackage{epsf}
\usepackage{epsfig}
\usepackage{setspace}
\usepackage{color}
\usepackage{algorithm}
\usepackage{algorithmic}
\usepackage{setspace}
\definecolor{SpringGreen4}{RGB}{0,139,69}

\begin{document}
\title{Multiple-GPU accelerated high-order gas-kinetic scheme for direct numerical simulation of compressible turbulence}

\author[BNU]{Yuhang Wang}
\ead{hskwyh@outlook.com}
\author[SUSTech1]{Guiyu Cao}
\ead{caogy@sustech.edu.cn}
\author[BNU]{Liang Pan\corref{cor}}
\ead{panliang@bnu.edu.cn}

\address[BNU]{Laboratory of Mathematics and Complex Systems, School of Mathematical Sciences, Beijing Normal University, Beijing, China}
\address[SUSTech1]{Academy for Advanced Interdisciplinary Studies, Southern University of Science and Technology, Shenzhen, China}
\cortext[cor]{Corresponding author}

\begin{abstract}
High-order gas-kinetic scheme (HGKS) has become a workable tool for
the direct numerical simulation (DNS) of turbulence. In this paper,
to accelerate the computation, HGKS is implemented with the
graphical processing unit (GPU) using the compute unified device
architecture (CUDA). Due to the limited available memory size, the
computational scale is constrained by single GPU. To conduct the
much large-scale DNS of turbulence, HGKS also be further upgraded
with multiple GPUs using message passing interface (MPI) and CUDA
architecture. The benchmark cases for compressible turbulence,
including Taylor-Green vortex and turbulent channel flows, are
presented to assess the numerical performance of HGKS with Nvidia
TITAN RTX and Tesla V100 GPUs. For single-GPU computation, compared
with the parallel central processing unit (CPU) code running on the
Intel Core i7-9700 with open multi-processing (OpenMP) directives,
7x speedup is achieved by TITAN RTX and 16x speedup is achieved by
Tesla V100. For multiple-GPU computation, multiple-GPU accelerated
HGKS code scales properly with the increasing number of GPU. The
computational time of parallel CPU code running on $1024$ Intel Xeon
E5-2692 cores with MPI is approximately $3$ times longer than that
of GPU code using $8$ Tesla V100 GPUs with MPI and CUDA. Numerical
results confirm the excellent performance of multiple-GPU
accelerated HGKS for large-scale DNS of turbulence. Alongside the
footprint in reduction of loading and writing pressure of GPU
memory, HGKS in GPU is also compiled with FP32 precision to evaluate
the effect of number formats precision. Reasonably, compared to the
computation with FP64 precision, the efficiency is improved and the
memory cost is reduced with FP32 precision. Meanwhile, the
differences in accuracy for statistical turbulent quantities appear.
For turbulent channel flows, difference in long-time statistical
turbulent quantities is acceptable between FP32 and FP64 precision
solutions. While the obvious discrepancy in instantaneous turbulent
quantities can be observed, which shows that FP32 precision is not
safe for DNS in compressible turbulence. The choice of precision
should depended on the requirement of accuracy and the available
computational resources.

\end{abstract}

\maketitle

\section{Introduction}
Turbulence is ubiquitous in natural phenomena and engineering
applications \cite{pope2001turbulent}. The understanding and
prediction of multiscale turbulent flow is one of the most difficult
problems for both mathematics and physical sciences. Direct
numerical simulation (DNS) solves the Navier-Stokes equations
directly, resolve all scales of the turbulent motion, and eliminate
modeling entirely \cite{kim1987turbulence}. With the advances of
high-order numerical methods and supercomputers, great success has
been achieved for the incompressible and compressible turbulent
flow, such as DNS in incompressible isotropic turbulence
\cite{chen1993far}, incompressible turbulent channel flow
\cite{kim1999turbulence, lee2015direct}, supersonic isotropic
turbulence \cite{wang2010hybrid,GKS-high-1}, compressible turbulent
channel flows \cite{coleman1995numerical, DNS-Li, yu2019genuine},
and  compressible flat plate turbulence from the supersonic to
hypersonic regime \cite{pirozzoli2004direct,lixinliangma8}.

In the past decades, the gas-kinetic scheme (GKS) has been developed
systematically based on the Bhatnagar-Gross-Krook (BGK) model
\cite{BGK-1,BGK-2} under the finite volume framework, and applied
successfully in the computations from low speed flow to hypersonic
one \cite{GKS-Xu1,GKS-Xu2}. The gas-kinetic scheme presents a gas
evolution process from kinetic scale to hydrodynamic scale, where
both inviscid and viscous fluxes are recovered from a time-dependent
and genuinely multi-dimensional gas distribution function at a cell
interface. Starting from a time-dependent flux function, based on
the two-stage fourth-order formulation \cite{GRP-high-1,GRP-high-2},
the high-order gas-kinetic scheme (HGKS) has been constructed and
applied for the compressible flow simulation
\cite{GKS-high-1,GKS-high-2,GKS-high-3,GKS-high-4}. The high-order
can be achieved with the implementation of the traditional
second-order GKS flux solver. More importantly, the high-order GKS
is as robust as the second-order scheme and works perfectly from the
subsonic to hypersonic viscous heat conducting flows. Originally,
the parallel HGKS code was developed with central processing unit
(CPU) using open multi-processing (OpenMP) directives. However, due
to the limited shared memory, the computational scale is constrained
for numerical simulation of turbulence. To perform the large-scale
DNS, the domain decomposition and the message passing interface
(MPI) \cite{MPI-1} are used for parallel implementation
\cite{GKS-high-DNS}. Due to the explicit formulation of HGKS, the
CPU code with MPI scales properly with the number of processors
used. The numerical results demonstrates the capability of HGKS as a
powerful DNS tool from the low speed to supersonic turbulence study
\cite{GKS-high-DNS}.

Graphical processing unit (GPU) is a form of hardware acceleration,
which is originally developed for graphics manipulation and is
extremely efficient at processing large amounts of data in parallel.
Since these units have a parallel computation capability inherently,
they can provide fast and low cost solutions to high performance
computing (HPC). In recent years, GPUs have gained significant
popularity in computational fluid dynamics \cite{GPU-1,GPU-2,GPU-3}
as a cheaper, more efficient, and more accessible alternative to
large-scale HPC systems with CPUs. Especially, to numerically study
the turbulent physics in much higher Reynolds number turbulence, the
extreme-scale DNS in turbulent flows has been implemented using
multiple GPUs, i.e., incompressible isotropic turbulence up to
$18432^3$ resolution \cite{yeung2020advancing} and turbulent pipe
flow up to $Re_{\tau} \approx 6000$ \cite{pirozzoli2021one}.
Recently, the three-dimension discontinuous Galerkin based HGKS has
been implemented in single-GPU computation using compute unified
device architecture (CUDA) \cite{GKS-GPU}. Obtained results are
compared with those obtained by Intel i7-9700 CPU using OpenMP
directives. The GPU code achieves 6x-7x speedup with TITAN RTX, and
10x-11x speedup with Tesla V100. The numerical results confirm the
potential of HGKS for large-scale DNS in turbulence.

A major limitation in single-GPU computation is its available
memory, which leads to a bottleneck in the maximum number of
computational mesh. In this paper, to implement much larger scale
DNS and accelerate the efficiency, HGKS is implemented with multiple
GPUs using CUDA and MPI architecture (MPI + CUDA). It is not
straightforward for the multiple-GPU programming, since the memory
is not shared across the GPUs and their tasks need to be coordinated
appropriately. The multiple GPUs are distributed across multiple
CPUs at the cost of having to coordinate GPU-GPU communication via
MPI. For WENO-based HGKS in single-GPU using CUDA, compared with the
CPU code using Intel Core i7-9700, 7x speedup is achieved for TITAN
RTX and 16x speedup is achieved for Tesla V100. In terms of the
computation with multiple GPUs, the HGKS code scales properly with
the increasing number of GPU. Numerical performance shows that the
data communication crossing GPU through MPI costs the relative
little time, while the computation time for flow field is the
dominant one in the HGKS code. For the MPI parallel computation, the
computational time of CPU code with supercomputer using 1024 Intel
Xeon E5-2692 cores is approximately $3$ times longer than that of
GPU code using $8$ Tesla V100 GPUs. It can be inferred that the
efficiency of GPU code with $8$ Tesla V100 GPUs approximately equals
to that of MPI code using $3000$ CPU cores in supercomputer. To
reduce the loading and writing pressure of GPU memory, the benefits
can be achieved by using FP32 (single) precision compared with FP64
(double) precision for memory-intensive computing
\cite{lehmann2021accuracy, haidar2020mixed}. Thence, HGKS in GPUs is
compiled with both FP32 precision and FP64 precision to evaluate the
effect of precision on DNS of compressible turbulence. As expect,
the efficiency can be improved and the memory cost can be reduced
with FP32 precision. However, the difference in accuracy between
FP32 and FP64 precision appears for the instantaneous statistical quantities of
turbulent channel flows, which is not negligible for long time simulations.

This paper is organized as follows. In Section 2, the high-order
gas-kinetic scheme is briefly reviewed. The GPU architecture and
code design are introduced in Section 3. Section 4 includes
numerical simulation and discussions. The last section is the
conclusion.

\section{High-order gas-kinetic scheme}
The three-dimensional BGK equation \cite{BGK-1,BGK-2} can be
written as
\begin{equation}\label{bgk}
f_t+uf_x+vf_y+wf_z=\frac{g-f}{\tau},
\end{equation}
where $\boldsymbol{u}=(u,v,w)^T$ is the particle velocity, $f$ is the
gas distribution function, $g$ is the three-dimensional Maxwellian
distribution and $\tau$ is the collision time. The collision term
satisfies the compatibility condition
\begin{equation} \label{compatibility}
\int \frac{g-f}{\tau}\psi \text{d}\Xi=0,
\end{equation}
where
$\displaystyle\psi=(1,u,v,w,\frac{1}{2}(u^2+v^2+w^2+\xi^2))^T$,
$\xi^2=\xi_1^2+...+\xi_N^2$,
$\text{d}\Xi=\text{d}u\text{d}v\text{d}w\text{d}\xi_1,...,\text{d}\xi_{N}$,
$N=(5-3\gamma)/(\gamma-1)$ is the internal degree of freedom,
$\gamma$ is the specific heat ratio and $\gamma=1.4$ is used in the
computation.

In this section, the finite volume scheme on orthogonal structured
mesh is provided as example. Taking moments of the BGK equation
Eq.\eqref{bgk} and integrating with respect to the cell
$\Omega_{ijk}$, the finite volume scheme can be expressed as
\begin{align}\label{semi}
\frac{\text{d}(Q_{ijk})}{\text{d}t}=\mathcal{L}(Q_{ijk}),
\end{align}
where $Q_{ijk}$ is the cell averaged conservative variable over
$\Omega_{ijk}$, and the operator $\mathcal{L}$ reads
\begin{equation}\label{semi2}
\mathcal{L}(Q_{ijk})=-\frac{1}{|\Omega_{ijk}|}\sum_{p=1}^6\mathbb{F}_{p}(t).
\end{equation}
where $\Omega_{ijk}$ is defined as
$\Omega_{ijk}=\overline{x}_i\times\overline{y}_j\times
\overline{z}_k$ with $\overline{x}_i=[x_i-\Delta x/2,x_i+\Delta
x/2], \overline{y}_j=[y_j-\Delta y/2,y_j+\Delta y/2],
\overline{z}_k=[z_k-\Delta z/2,z_k+\Delta z/2]$, and
$\mathbb{F}_{p}(t)$ is the numerical flux across the cell interface
$\Sigma_{p}$. The numerical flux in $x$-direction is given as
example
\begin{align*}
\mathbb{F}_{p}(t)=\iint_{\Sigma_{p}}
F(Q)\cdot\boldsymbol{n}\text{d}\sigma=\sum_{m,n=1}^2\omega_{mn}
\int\psi u
f(\boldsymbol{x}_{i+1/2,j_m,k_n},t,\boldsymbol{u},\xi)\text{d}\Xi\Delta y\Delta z,
\end{align*}
where $\boldsymbol{n}$ is the outer normal direction. The Gaussian
quadrature is used over the cell interface, where $\omega_{mn}$ is
the quadrature weight,
$\boldsymbol{x}_{i+1/2,m,n}=(x_{i+1/2},y_{j_m},z_{k_n})$ and
$(y_{j_m},z_{k_n})$ is the Gauss quadrature point of cell interface
$\overline{y}_j\times\overline{z}_k$. Based on the integral solution
of BGK equation Eq.\eqref{bgk}, the gas distribution function
$f(\boldsymbol{x}_{i+1/2,j_m,k_n},t,\boldsymbol{u},\xi)$ in the
local coordinate can be given by
\begin{equation}
f(\boldsymbol{x}_{i+1/2,j_m,k_n},t,\boldsymbol{u},\xi)=\frac{1}{\tau}\int_0^t
g(\boldsymbol{x}',t',\boldsymbol{u}, \xi)e^{-(t-t')/\tau}\text{d}t'+e^{-t/\tau}f_0(-\boldsymbol{u}t,\xi),
\end{equation}
where
$\boldsymbol{x}'=\boldsymbol{x}_{i+1/2,j_m,k_n}-\boldsymbol{u}(t-t')$
is the trajectory of particles, $f_0$ is the initial gas
distribution function, and $g$ is the corresponding equilibrium
state. With the first order spatial derivatives, the second-order
gas distribution function at cell interface can be expressed as
\begin{align}\label{flux}
f(\boldsymbol{x}_{i+1/2,j_m,k_n},t,\boldsymbol{u},\xi)=&(1-e^{-t/\tau})g_0+
((t+\tau)e^{-t/\tau}-\tau)(\overline{a}_1u+\overline{a}_2v+\overline{a}_3w)g_0\nonumber\\
+&(t-\tau+\tau e^{-t/\tau}){\bar{A}} g_0\nonumber\\
+&e^{-t/\tau}g_r[1-(\tau+t)(a_{1}^{r}u+a_{2}^{r}v+a_{3}^{r}w)-\tau A^r)](1-H(u))\nonumber\\
+&e^{-t/\tau}g_l[1-(\tau+t)(a_{1}^{l}u+a_{2}^{l}v+a_{3}^{l}w)-\tau A^l)]H(u),
\end{align}
where the equilibrium state $g_{0}$ and the corresponding
conservative variables $Q_{0}$ can be determined by the
compatibility condition
\begin{align*}
\int\psi g_{0}\text{d}\Xi=Q_0=\int_{u>0}\psi
g_{l}\text{d}\Xi+\int_{u<0}\psi g_{r}\text{d}\Xi.
\end{align*}
The following numerical tests on the compressible turbulent flows
without discontinuities will be presented, thus the collision time
for the flow without discontinuities takes
\begin{align*}
\tau=\frac{\mu}{p},
\end{align*}
where $\mu$ is viscous coefficient and $p$ is the pressure at cell
interface determined by $Q_0$. With the reconstruction of
macroscopic variables, the coefficients in Eq.\eqref{flux} can be
fully determined by the reconstructed derivatives and compatibility
condition
\begin{equation*}
    \begin{aligned}
\displaystyle \langle a_{1}^{k}\rangle=\frac{\partial
Q_{k}}{\partial x}, \langle
a_{2}^{k}\rangle=\frac{\partial Q_{k}}{\partial y},
\langle a_{3}^{k}\rangle&=\frac{\partial Q_{k}}{\partial
z}, \langle
a_{1}^{k}u+a_{2}^{k}v+a_{3}^{k}w+A^{k}\rangle=0,\\ \displaystyle
\langle\overline{a}_1\rangle=\frac{\partial Q_{0}}{\partial
x}, \langle\overline{a}_2\rangle=\frac{\partial
Q_{0}}{\partial y},
\langle\overline{a}_3\rangle&=\frac{\partial Q_{0}}{\partial
z},
\langle\overline{a}_1u+\overline{a}_2v+\overline{a}_3w+\overline{A}\rangle=0,
\end{aligned}
\end{equation*}
where $k=l,r$ and $\langle...\rangle$ are the moments of the
equilibrium $g$ and defined by
\begin{align*}
\langle...\rangle=\int g (...)\psi \text{d}\Xi.
\end{align*}
More details of the gas-kinetic scheme can be found in the
literature \cite{GKS-Xu1,GKS-Xu2}. Thus, the numerical flux can be
obtained by taking moments of the gas distribution function
Eq.\eqref{flux}, and the semi-discretized finite volume scheme
Eq.\eqref{semi} can be fully given.

Recently, based on the time-dependent flux function of the
generalized Riemann problem solver (GRP)
\cite{GRP-high-1,GRP-high-2} and gas-kinetic scheme
\cite{GKS-high-2,GKS-high-3,GKS-high-4}, a two-stage fourth-order
time-accurate discretization was recently developed for Lax-Wendroff
type flow solvers. A reliable framework was provided to construct
higher-order gas-kinetic scheme, and the high-order scheme is as
robust as the second-order one and works perfectly from the subsonic
to hypersonic flows. In this study, the two-stage method is used for
temporal accuracy. Consider the following time dependent equation
\begin{align*}
\frac{\partial Q}{\partial t}=\mathcal {L}(Q),
\end{align*}
with the initial condition at $t_n$, i.e., $Q(t=t_n)=Q^n$,
where $\mathcal {L}$ is an operator for spatial derivative of flux,
the state $Q^{n+1}$ at $t_{n+1}=t_n+\Delta t$  can be updated with
the following formula
\begin{equation}\label{two-stage}
\begin{split}
&Q^*=Q^n+\frac{1}{2}\Delta t\mathcal {L}(Q^n)+\frac{1}{8}\Delta
t^2\partial_t\mathcal{L}(Q^n), \\
Q^{n+1}=&Q^n+\Delta t\mathcal {L}(Q^n)+\frac{1}{6}\Delta
t^2\big(\partial_t\mathcal{L}(Q^n)+2\partial_t\mathcal{L}(Q^*)\big).
\end{split}
\end{equation}
It can be proved that for hyperbolic equations the above temporal
discretization Eq.\eqref{two-stage} provides a fourth-order time accurate solution for
$Q^{n+1}$. To implement two-stage fourth-order method for
Eq.\eqref{semi}, a linear function is used to approximate the time
dependent numerical flux
\begin{align}\label{expansion-1}
\mathbb{F}_{p}(t)\approx\mathbb{F}_{p}^n+ \partial_t
\mathbb{F}_{p}^n(t-t_n).
\end{align}
Integrating Eq.\eqref{expansion-1} over $[t_n, t_n+\Delta t/2]$ and
$[t_n, t_n+\Delta t]$, we have the following two equations
\begin{equation*}
\begin{aligned}
\mathbb{F}_{p}^n\Delta t&+\frac{1}{2}\partial_t
\mathbb{F}_{p}^n\Delta t^2 =\int_{t_n}^{t_n+\Delta t}\mathbb{F}_{p}(t)\text{d}t, \\
\frac{1}{2}\mathbb{F}_{p}^n\Delta t&+\frac{1}{8}\partial_t
\mathbb{F}_{p}^n\Delta t^2 =\int_{t_n}^{t_n+\Delta t/2}\mathbb{F}_{p}(t)\text{d}t.
\end{aligned}
\end{equation*}
The coefficients $\mathbb{F}_{p}^n$ and $\partial_t
\mathbb{F}_{p}^n$ at the initial stage can be determined by solving
the linear system. According to Eq.\eqref{semi2}, $\mathcal
{L}(Q_{i}^n)$ and the temporal derivative $\partial_t\mathcal
{L}(Q_{i}^n)$ at $t^n$ can be constructed by
\begin{align*}
\mathcal{L}(Q_{i}^n)&=-\frac{1}{|\Omega_{i}|}\sum_{p=1}^6\mathbb{F}_{p}^n,~~
\partial_t\mathcal{L}(Q_{i}^n)=-\frac{1}{|\Omega_{i}|}\sum_{p=1}^6\partial_t\mathbb{F}_{p}^n.
\end{align*}
The flow variables $Q^*$ at the intermediate stage can be updated.
Similarly, $\mathcal{L}(Q_{i}^*)$, $\partial_t\mathcal {L}(Q_{i}^*)$
at the intermediate state can be constructed and $Q^{n+1}$ can be
updated as well. For the high-order spatial accuracy, the
fifth-order WENO method \cite{WENO-JS, WENO-Z} is adopted. For the
three-dimensional computation, the dimension-by-dimension
reconstruction is used for HGKS \cite{GKS-high-3}.

\section{HGKS code design on GPU}
CPUs and GPUs are equipped with different architectures and are
built for different purposes. The CPU is suited to general
workloads, especially those for which per-core performance are
important. CPU is designed to execute complex logical tasks quickly,
but is limited in the number of threads. Meanwhile, GPU is a form of
hardware acceleration, which is originally developed for graphics
manipulation and is extremely efficient at processing large amounts
of data in parallel. Currently, GPU has gained significant
popularity in high performance scientific computing
\cite{GPU-1,GPU-2,GPU-3}. In this paper, to accelerate the
computation, the WENO based HGKS will be implemented with single GPU
using CUDA. To conduct the large-scale DNS in turbulence
efficiently, the HGKS also be further upgraded with multiple GPUs using MPI + CUDA.

\subsection{Single-GPU accelerated HGKS}
\begin{figure}[!h]
    \centering
    \includegraphics[width=0.85\textwidth]{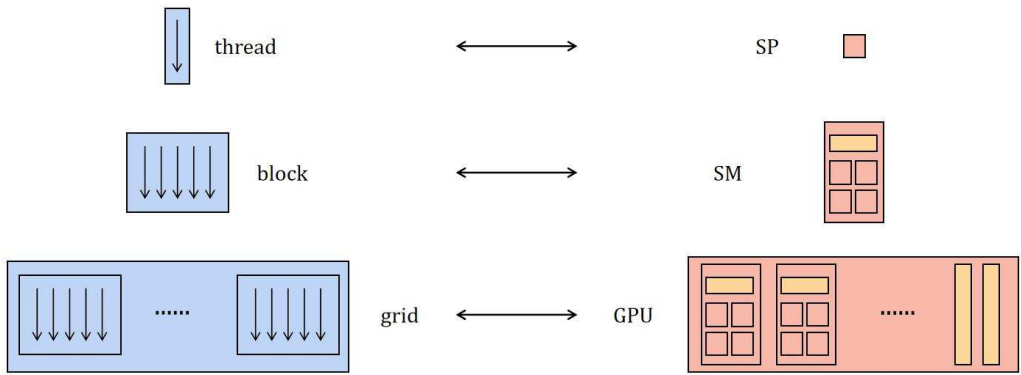}
    \caption{\label{cuda-gpu} Correspondence between software CUDA and hardware GPU.}
\end{figure}

The CPU is regarded as host, and GPU is treated as device.
Data-parallel, compute-intensive operations running on the host are
transferred to device by using kernels, and kernels are executed on
the device by many different threads. For CUDA, these threads are
organized into thread blocks, and thread blocks constitute a grid.
Such computational structures build connection with Nvidia GPU
hardware architecture. The Nvidia GPU consists of multiple streaming
multiprocessors (SMs), and each SM contains streaming processors
(SPs). When invoking a kernel, the blocks of grid are distributed to
SMs, and the threads of each block are executed by SPs. In summary,
the correspondence between software CUDA and hardware Nvidia GPU is
simply shown in Fig.\ref{cuda-gpu}. In the following, the single-GPU
accelerated HGKS will be introduced briefly.

\begin{algorithm}[!h]
    \setstretch{1.25}
    \begin{algorithmic}
        \STATE
        \textbf{\color{blue} Initialization}
        \WHILE {TIME $ \leq $ TSTOP}
        \STATE dimGrid $=$ dim3($N_z$, $N_x/\text{block}_x$, 1)
        \STATE dimBlock $=$ dim3($\text{block}_x$,$N_y$,1)
        \STATE \textbf{STEP 1} : \textbf{\color{blue} Calculation of time step}
        \STATE \textbf{CALL} GETTIMESTEP$<<<$dimGrid,dimBlock$>>>$
        \STATE istat=cudaDeviceSynchronize()
        \STATE \textbf{STEP 2} : \textbf{\color{blue} WENO reconstruction}
        \STATE \textbf{CALL} WENO-x$<<<$dimGrid,dimBlock$>>>$
        \STATE istat=cudaDeviceSynchronize()
        \STATE \textbf{STEP 3} : \textbf{\color{blue} Computation of flux}
        \STATE \textbf{CALL} FLUX-x$<<<$dimGrid,dimBlock$>>>$   \\
        \STATE istat=cudaDeviceSynchronize() \\
        \quad  \textcolor{SpringGreen4}{\%  kernel for WENO reconstruction and flux calculation in $x$ direction as example.}   \\
        \quad  \textcolor{SpringGreen4}{\%  reconstruction and flux calculation in $y$ and $z$ directions can be also implemented.}
        \STATE \textbf{STEP 4} : \textbf{\color{blue} Update of flow variables}
        \ENDWHILE \\
    \end{algorithmic}
    \caption{\label{GPU-algorithm} Single-GPU HGKS code using CUDA}
\end{algorithm}
\begin{figure}[!htp]
    \centering
    \includegraphics[width=0.55\textwidth]{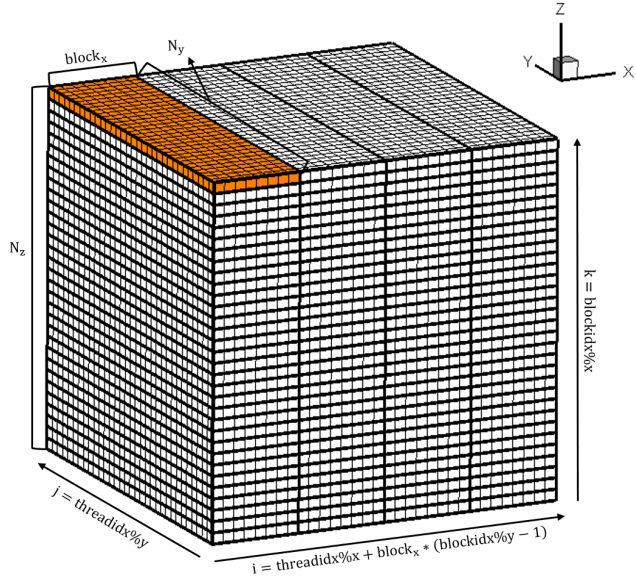}
    \caption{\label{thread-3d} Specifications of grid for the three-dimensional structured meshes.}
\end{figure}

Device implementation is handled automatically by the CUDA, while
the kernels and grids need to be specified for HGKS. One updating
step of single-GPU accelerated HGKS code can be divided into several
kernels, namely, initialization, calculation of time step, WENO
reconstruction, flux computation at cell interface in $x$, $y$, $z$,
and the update of flow variables. The main parts of the single-GPU
accelerated HGKS code are labeled in blue as
Algorithm.\ref{GPU-algorithm}. Due to the identical processes of
calculation for each kernel, it is natural to use one thread for one
computational cell when setting grids. Due to the limited number of
threads can be obtained in one block (i.e., the maximum $1024$ threads in one block), the whole computational
domain is required to be divided into several parts. Assume that the total number of cells for
computational mesh is $N_x\times N_y\times N_z$, which is
illustrated in Fig.\ref{thread-3d}. The computational domain is
divided into $N_x/\text{block}_x$ parts in $x$-direction, where the choice of $\text{block}_x$ for each kernel is an integer defined according to the requirement and experience. If $N_x$ is not divisible by $\text{block}_x$, an extra block is needed. When
setting grid, the variables "dimGrid" and "dimBlock" are defined as
\begin{align*}
{\rm dimGrid}&={\rm dim3}(N_z, N_x/\text{block}_x, 1), \\
 {\rm dimBlock}&={\rm dim3}(\text{block}_x,N_y,1).
\end{align*}
As shown in Fig.\ref{thread-3d}, one thread is assigned to complete
the computations of a cell $\Omega_{ijk}$, and the one-to-one
correspondence of thread index $\text{threadidx}$, block index  $\text{blockidx}$ and cell index $(i,j,k)$ is given by
\begin{align*}
i&=\text{threadidx}\%\text{x}+\text{block}_x*(\text{blockidx}\%\text{y}-1),\\
j&=\text{threadidx}\%\text{y}, \\
k&= \text{blockidx}\%\text{x}.
\end{align*}
Thus, the single-GPU accelerated HGKS code can be implemented after
specifying kernels and grids.

\begin{figure}[!h]
\centering
\includegraphics[width=0.9\textwidth]{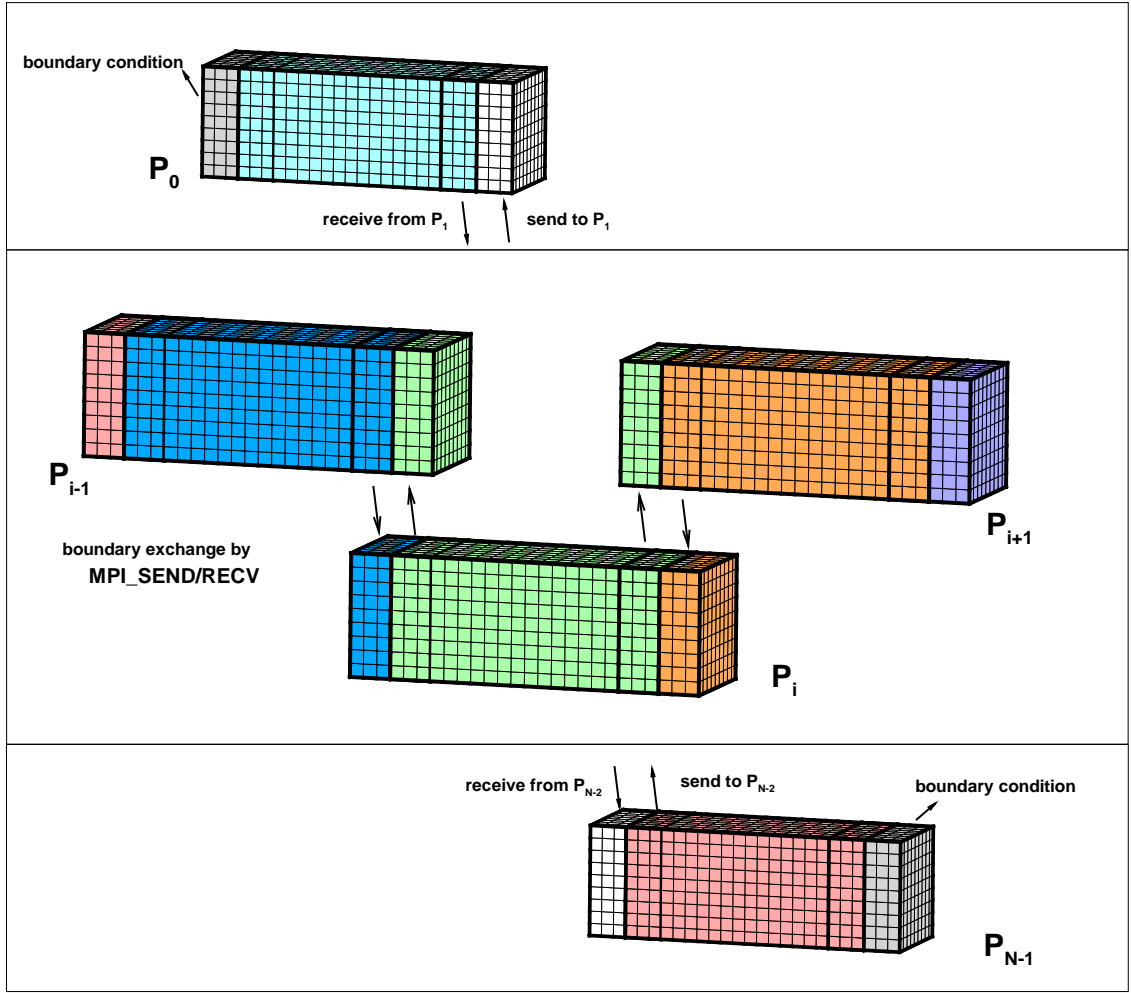}
\caption{\label{total-separate} The domain decomposition and GPU-GPU communications for multiple-GPU implementation.}
\end{figure}

\subsection{Multiple-GPU accelerated HGKS}
Due to the limited computational resources of single-GPU, it is
natural to develop the multiple-GPU accelerated HGKS for large-scale
DNS of turbulence. The larger computational scale can be achieved
with the increased available device memories. It is not
straightforward for the programming with multiple-GPUs, because the
device memories are not shared across GPUs and the tasks need to be
coordinated appropriately. The computation with multiple GPUs is
implemented using multiple CPUs and multiple GPUs, where the GPU-GPU
communication is coordinated via the MPI. CUDA-aware MPI library is
chosen \cite{CUDA-MPI}, where GPU data can be directly passed by MPI
function and transferred in the efficient way automatically.

\begin{algorithm}[!h]
\setstretch{1.25}
\begin{algorithmic}
\IF{($i =0$ ~or  ~$i = N-1$)}
\STATE \textbf{\color{blue} Boundary conditions} for $P_0$ and $P_{N-1}$\\
\quad  \textcolor{SpringGreen4}{\%  for wall boundary, boundary condition can be given directly.}  \\
\quad  \textcolor{SpringGreen4}{\%  for periodic boundary, \textbf{\color{blue}MPI$\_$SEND}  and  \textbf{\color{blue}MPI$\_$RECV} are also needed.}
\ENDIF \\
\IF{($i =$ even number)}
\STATE \textbf{\color{blue}MPI$\_$SEND} data from $P_{i}$ to $P_{i-1}$,
\STATE \textbf{\color{blue}MPI$\_$RECV} data from $P_{i-1}$ to $P_{i}$.
\ELSE
\STATE \textbf{\color{blue}MPI$\_$RECV} data from $P_{i+1}$  to $P_{i}$,
\STATE \textbf{\color{blue}MPI$\_$SEND} data from $P_{i}$ to $P_{i+1}$
\ENDIF
\IF{($i =$ odd number)}
\STATE \textbf{\color{blue}MPI$\_$SEND} data from $P_{i}$ to $P_{i-1}$,
\STATE \textbf{\color{blue}MPI$\_$RECV} data from $P_{i-1}$ to $P_{i}$.
\ELSE
\STATE \textbf{\color{blue}MPI$\_$RECV} data from $P_{i+1}$  to $P_{i}$,
\STATE \textbf{\color{blue}MPI$\_$SEND} data from $P_{i}$ to $P_{i+1}$
\ENDIF
\end{algorithmic}
\caption{\label{Boundary-algorithm} Data communication for MPI+CUDA}
\end{algorithm}

\begin{figure}[htp]
\centering
\includegraphics[width=0.9\textwidth]{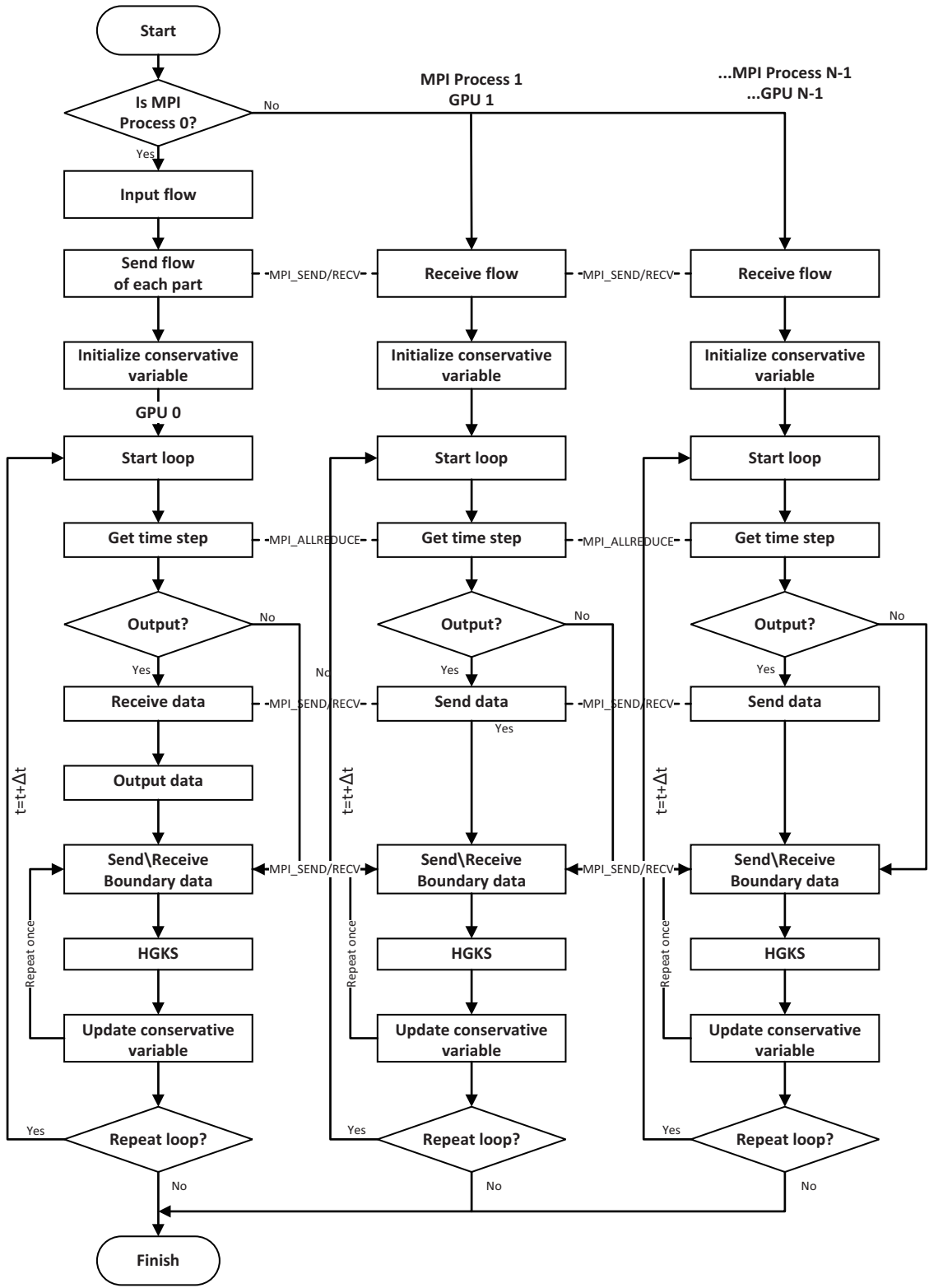}
\caption{\label{code-frame} Multiple-GPU accelerated HGKS code using MPI+CUDA.}
\end{figure}

For the parallel computation with MPI, the domain decomposition are
required for both CPU and GPU, where both one-dimensional (slab) and
two-dimensional (pencil) decomposition are allowed. For CPU
computation with MPI, hundreds or thousands of CPU cores are usually
used in the computation, and the pencil decomposition can be used to
improve the efficiency of data communications \cite{GKS-high-DNS}. In terms of GPU
computation with MPI, the slab decomposition is used, which reduces
the frequency of GPU-GPU communications and improves performances.
The domain decomposition and GPU-GPU communications are presented in
Fig.\ref{total-separate}, where the computational domain is divided
into $N$ parts in $x$-direction and $N$ GPUs are used. For
simplicity, $P_i$ denotes the $i$-th MPI process, which deals with
the tasks of $i$-th decomposed computational domain in the $i$-th
GPU. To implement the WENO reconstruction, the processor $P_i$
exchanges the data with $P_{i-1}$ and $P_{i+1}$ using MPI$\_$RECV
and MPI$\_$SEND. The wall boundary condition of $P_0$ and $P_{N-1}$
can be implemented directly. For the periodic boundary, the data is
exchanged between $P_0$ and $P_{N-1}$ using MPI$\_$RECV and
MPI$\_$SEND as well. In summary, the detailed data communication for
GPU computation with MPI are shown in
Algorithm.\ref{Boundary-algorithm}. Fig.\ref{code-frame} shows the
HGKS code frame with multiple-GPUs. The initial data will be divided
into $N$ parts according to domain decomposition and sent to
corresponding process $P_1$ to $P_{N-1}$ using MPI$\_$SEND. As for
input and output, the MPI process $P_0$ is responsible for
distributing data to other processors and collecting data from other
processors, using MPI$\_$SEND and MPI$\_$RECV. For each decomposed
domain, the computation is executed with the prescribed GPU. The
identical code is running for each GPU, where kernels and grids are
similarly set as single-GPU code \cite{GKS-GPU}.

\subsection{Detailed parameters of CPUs and GPUs}
In following studies, the CPU codes run on the desktop with Intel
Core i7-9700 CPU using OpenMP directives and TianHe-II supercomputer
system with MPI and Intel mpiifort compiler, The HGKS codes with CPU
are all compiled with FP64 precision. For the GPU computations, the
Nvidia TITAN RTX GPU and Nvidia Tesla V100 GPU are used with Nvidia
CUDA and Nvidia HPC SDK. The detailed parameters of CPU and GPU
are given in Table.\ref{GPU-CPU-TG-A} and Table.\ref{GPU-CPU-TG-B},
respectively.  For TITAN RTX GPU, the GPU-GPU communication is
achieved by connection traversing PCIe, and there are $2$ TITAN TRTX
GPUs in one node. For Tesla V100 GPU, there are $8$ GPUs inside one
GPU node, and more nodes are needed for more than $8$ GPUs. The
GPU-GPU communication in one GPU node is achieved by Nvidia NVLink.
The communication across GPU nodes can be achieved by Remote Direct
Memory Access (RDMA) technique, including iWARP, RDMA over Converged
Ethernet (RoCE) and InfiniBand. In this paper, RoCE is used for
communication across GPU nodes.

As shown in Table.\ref{GPU-CPU-TG-B}, Tesla V100 is equipped with
more FP64 precision cores and much stronger FP64 precision
performance than TITAN RTX. For FP32 precision performance, two GPUs
are comparable. Nvidia TITAN RTX is powered by Turing architecture,
while Nvidia Tesla V100 is built to HPC by advanced Ampere
architecture. Additionally, Tesla V100 outweighs TITAN RTX in GPU
memory size and memory bandwidth. Thence,  much excellent
performance is excepted in FP64 precision simulation with Tesla
V100. To validate the performance with FP32 and FP64 precision,
numerical examples will be presented.

\begin{table}[!h]
\centering
\begin{tabular}{c|c|c|c}
\hline
\hline
~             &    CPU       &  Clock rate      & Memory size   \\
\hline
Desktop   &    Intel Core i7-9700      &$3.0$ GHz &$16$ GB/$8$ cores \\
\hline
TianHe-II &   Intel Xeon E5-2692 v2     &$2.2$ GHz   &$64$ GB/$24$ cores\\
\hline
\hline
\end{tabular}
\caption{\label{GPU-CPU-TG-A} The detailed parameters of CPUs.}
~\\~\\
\centering
\begin{tabular}{c|c|c}
\hline
\hline
~                                     &  Nvidia TITAN RTX & Nvidia Tesla V100  \\
\hline
Clock rate                       &  1.77 GHz                &1.53 GHz \\
\hline
Stream multiprocessor    &  72                            &80  \\
\hline
FP64 core per SM         &   2                             &32 \\
\hline
FP32 precision performance           &   16.3 Tflops           &15.7 Tflops\\
\hline
FP64 precision performance            &   509.8 Gflops           &7834 Gflops \\
\hline
GPU memory size            &   24 GB           &32 GB \\
\hline
Memory bandwidth            &   672 GB/s           &897 GB/s \\
\hline
\hline
\end{tabular}
\caption{\label{GPU-CPU-TG-B}  The detailed parameters of GPUs.}
\end{table}

\section{Numerical tests and discussions}
Benchmarks for compressible turbulent flows, including Taylor-Green
vortex and turbulent channel flows, are presented to validate the
performance of multiple-GPU accelerated HGKS. In this section, we
mainly concentrate on the computation with single-GPU and
multiple-GPUs. The GPU-CPU comparison in terms of computational
efficiency will be presented firstly, and the scalability of
parallel multiple-GPU code will also be studied. Subsequently, HGKS
implemented with GPUs is compiled with both FP32 and FP64 precision
to evaluate the effect of precision on DNS of compressible
turbulence. The detailed comparison on numerical efficiency and
accuracy of statistical turbulent quantities with FP32 precision and
FP64 precision is also presented.

\subsection{Taylor-Green vortex for GPU-CPU comparison}
Taylor-Green vortex (TGV) is a classical problem in fluid dynamics
developed to study vortex dynamics, turbulent decay and energy
dissipation process \cite{Case-Brachet,Case-Debonis,Case-Bull}. In
this case, the detailed efficiency comparisons of HGKS running on
the CPU and GPU will be given. The flow is computed within a
periodic square box defined as $-\pi L\leq x, y, z\leq \pi L$. With
a uniform temperature, the initial velocity and pressure are given
by
\begin{equation*}
    \begin{aligned}
        U=&V_0\sin(\frac{x}{L})\cos(\frac{y}{L})\cos(\frac{z}{L}),\\
        V=&-V_0\cos(\frac{x}{L})\sin(\frac{y}{L})\cos(\frac{z}{L}),\\
        W=&0,\\
        p=&p_0+\frac{\rho_0V_0^2}{16}(\cos(\frac{2x}{L})+\cos(\frac{2y}{L}))(\cos(\frac{2z}{L})+2).
    \end{aligned}
\end{equation*}
In the computation, $L=1, V_0=1, \rho_0=1$, and the Mach number
takes $Ma = V_0/c_0=0.1$, where $c_0$ is the sound speed. The fluid
is a perfect gas with $\gamma=1.4$, Prandtl number $Pr=1.0$, and
Reynolds number $Re=1600$. The characteristic convective time $t_c =
L/V_0$.

\begin{table}[!h]
    \centering
    \begin{tabular}{c|c|c|c|c|c|c}
        \hline
        \hline
        Mesh size    & Time step     &  i7-9700  & TITAN RTX & $S_{gc}$   & Tesla V100 & $S_{gc}$\\
        \hline
        $64^3$  & $3.571\times10^{-3}$   &0.757    &0.126      &6.0 &0.047& 16.1  \\
        \hline
        $128^3$  & $1.785\times10^{-3}$  & 11.421   &1.496  &7.6 &0.714& 16.0\\
        \hline
        $256^3$  & $8.925\times10^{-4}$ &182.720   &24.256  &7.4 &11.697&15.4\\
        \hline
        \hline
    \end{tabular}
    \caption{\label{GPU-CPU-TG-C} Taylor-Green vortex: the total execution  times and speedup for
single-GPU versus CPU.}
\end{table}

In this section, GPU codes are all compiled with FP64 precision. To
test the performance of single-GPU code, this case is run with
Nvidia TITAN RTX and Nvidia Tesla V100 GPUs using Nvidia CUDA. As
comparison, the CPU code running on the Intel Core i7-9700 is also
tested. The execution times with different meshes are shown in
Table.\ref{GPU-CPU-TG-C}, where the total execution time for CPU and
GPUs are given in terms of hours and $20$ characteristic convective
time are simulated. Due to the limitation of single GPU memory size,
the uniform meshes with $64^3$, $128^3$ and $256^3$ cells are used.
For the CPU computation, $8$ cores are utilized with OpenMP parallel
computation, and the corresponding execution time
$T_{8,CPU}^{total}$ is used as the base for following comparisons.
The speedups $S_{gc}$ are also given in Table.\ref{GPU-CPU-TG-C},
which is defined as
\begin{equation*}
    \displaystyle {S_{gc}}=\frac{T_{8, CPU}^{total}}{T_{1, GPU}^{total}},
\end{equation*}
where $T_{8, CPU}^{total}$ and $T_{1, GPU}^{total}$ denote the total
execution time of $8$ CPU cores and $1$ GPU, respectively. Compared
with the CPU code, 7x speedup is achieved by single TITAN RTX GPU
and 16x speedup is achieved by single Tesla V100 GPU. Even though
the Tesla V100 is approximately $15$ times faster than the TITAN RTX
in FP64 precision computation ability, Tesla V100 only performs $2$
times faster than TITAN RTX. It indicates the memory bandwidth is
the bottleneck for current memory-intensive simulations, and the
detailed technique for GPU programming still required to be
investigated.

\begin{table}[!h]
\centering
\begin{tabular}{c|c|c|c}
\hline
\hline
 Mesh size & Time step            &No.  CPUs & $T_{n,CPU}^{total}$    \\
\hline
 $128^3$  & $1.785\times10^{-3}$  & 16    & 13.3   \\
\hline
$256^3$  & $8.925\times10^{-4}$  & 256   & 13.5  \\
\hline
$512^3$  &$4.462\times10^{-4}$   & 1024  & 66    \\
\hline
\hline
\end{tabular}
\caption{\label{time_table2} Taylor-Green vortex: the detailed
computational parameters in TianHe-II supercomputers.}
\end{table}

\begin{table}[!h]
\centering
\begin{tabular}{c|c|c|c|c}
\hline
\hline
Mesh  size & No. GPUs    & $T_{n,GPU}^{total}$  & $T_{n,GPU}^{flow}$   & $T_{n,GPU}^{com}$  \\
\hline
$256^3$   & 1    &11.697  & ~  & ~  \\
\hline
$256^3$   & 2    & 6.201 &6.122  & 0.079 \\
\hline
$256^3$   & 4    & 3.223  &2.997  & 0.226 \\
\hline
$256^3$   & 8    & 1.987 &1.548  & 0.439 \\
\hline
\hline
$512^3$  & 8    & 22.252  &21.383  &  0.869 \\
\hline
$512^3$  & 12    & 19.055 & 17.250 &  1.805  \\
\hline
$512^3$  & 16    & 13.466 &  11.483 &  1.983 \\
\hline
\hline
\end{tabular}
\caption{\label{GPU-CPU-TG-D} Taylor-Green vortex: the detailed
computational parameters with Tesla V100.}
\end{table}

\begin{figure}[!h]
\centering
\includegraphics[width=0.55\textwidth]{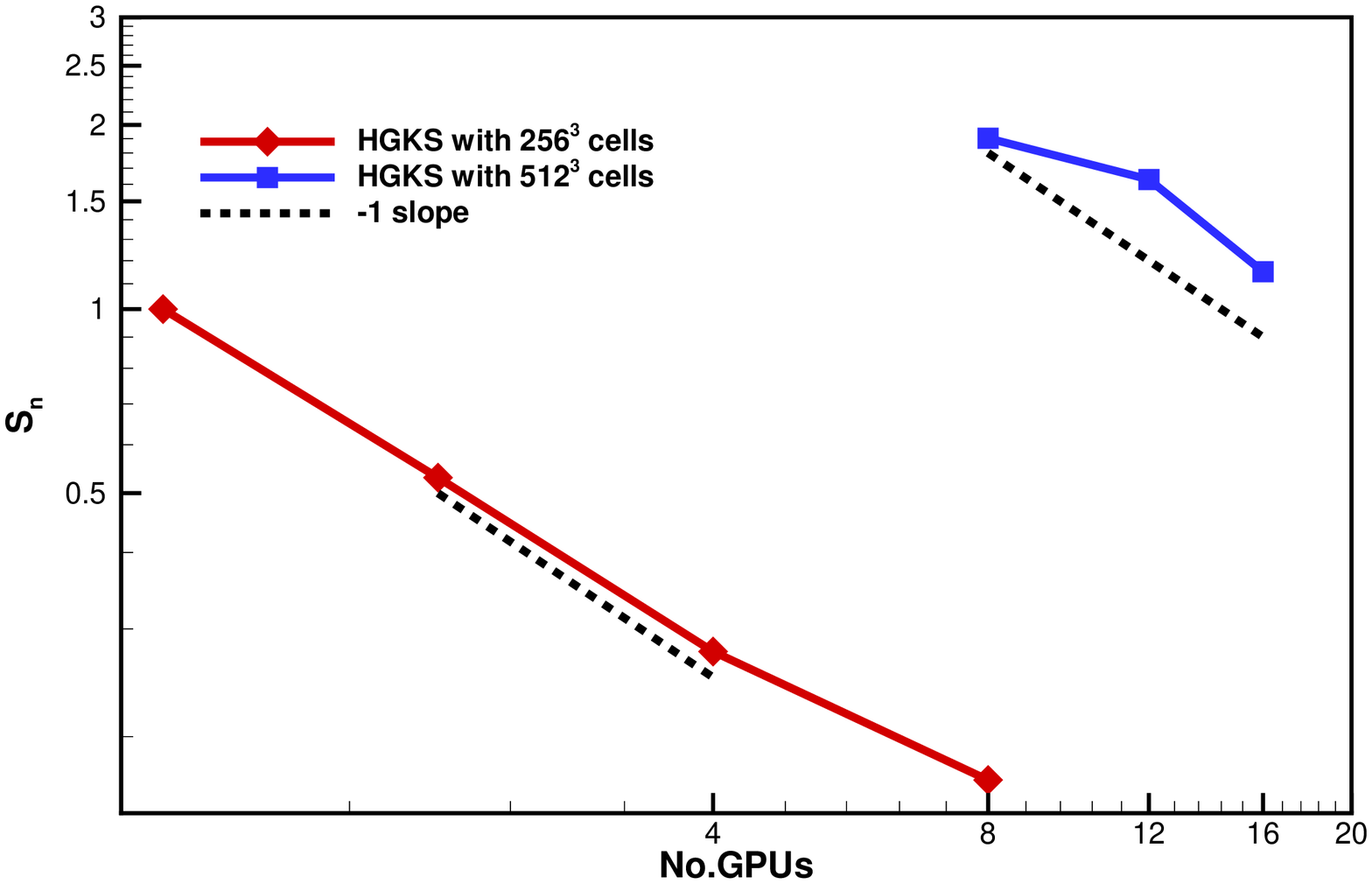}
\caption{\label{mpi-efficience} Taylor-Green vortex: speedup and
efficiency of Tesla V100 GPUs with $256^3$ and $512^3$ uniform
cells.}
~\\
\includegraphics[width=0.495\textwidth]{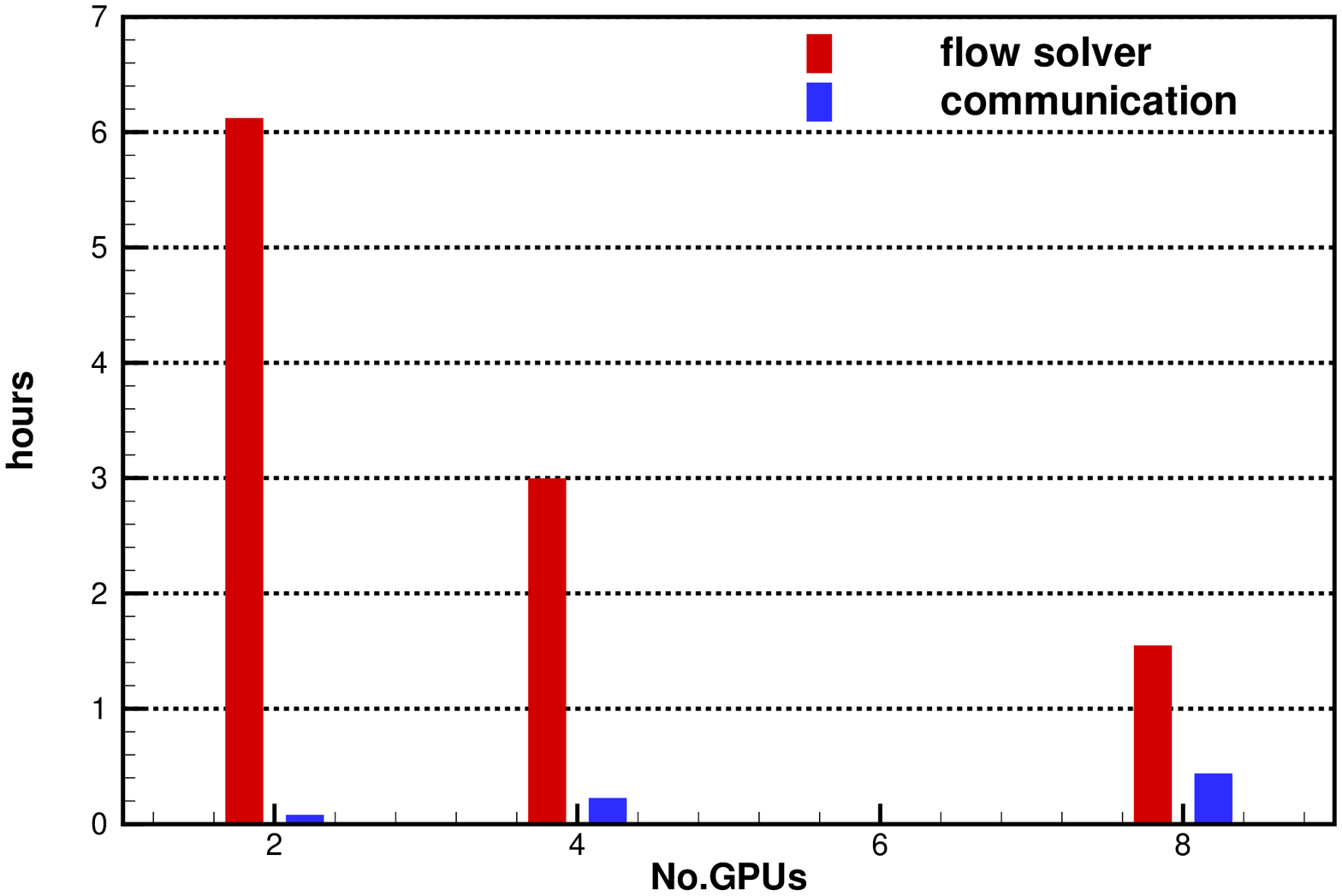}
\includegraphics[width=0.495\textwidth]{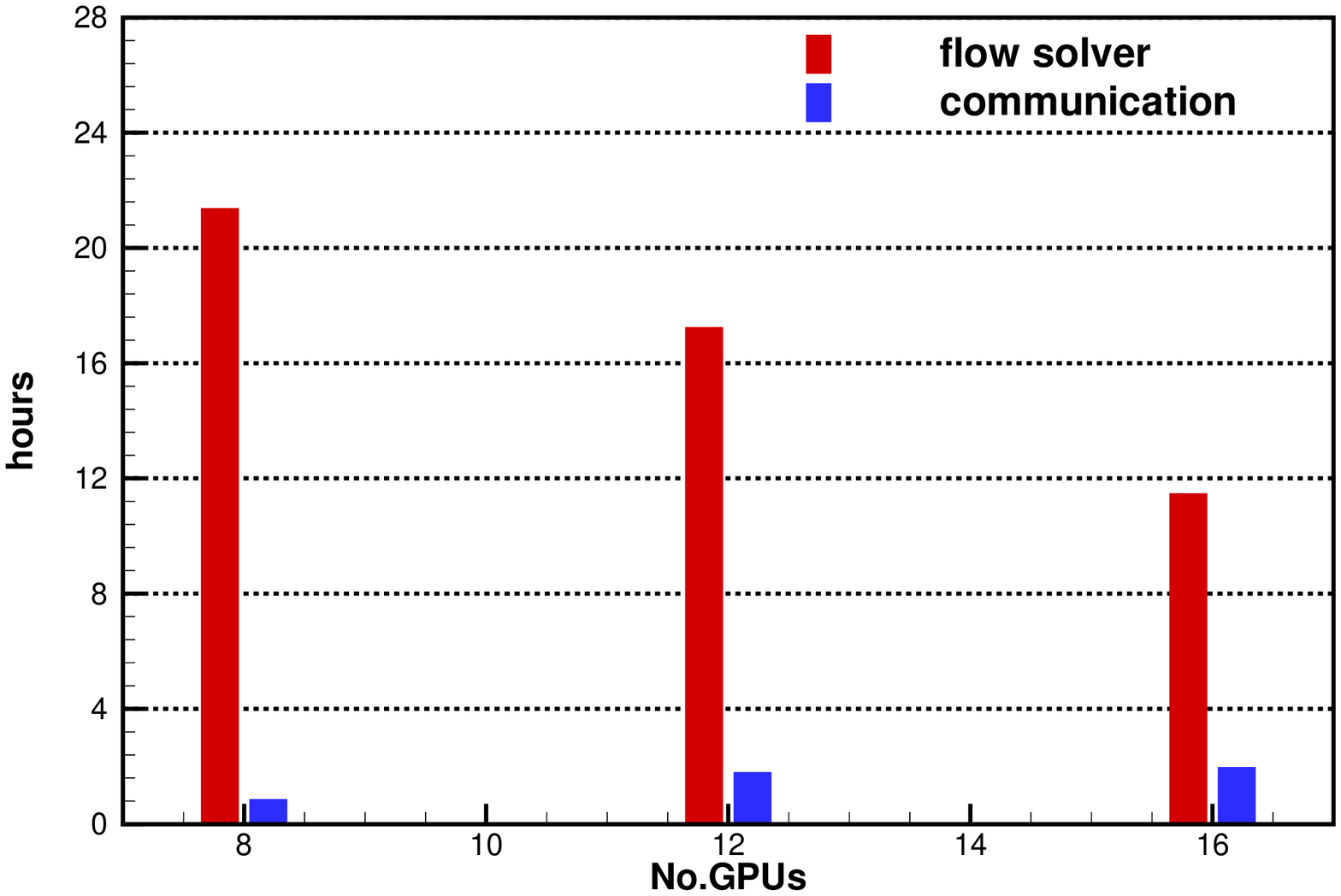}
\caption{\label{mpi-efficience-2} Taylor-Green vortex: the histogram
of $T_{com}^n$ and $T_{flow}^n$ using different number of Tesla V100
GPUs with $256^3$ (left) and $512^3$ (right) uniform cells.}
\end{figure}

For single TITAN TRX with $24$ GB memory size, the maximum
computation scale is $256^3$ cells. To enlarge the computational
scale, the multiple-GPU accelerated HGKS are designed. Meanwhile,
the parallel CPU code has been run on the TianHe-II supercomputer
\cite{GKS-high-DNS}, and the execution times with respect to the
number of CPU cores are presented  in Table.\ref{time_table2} for
comparison.  For the case with $256^3$ cells, the computational time
of single Tesla V100 GPU is comparable with the MPI code with
supercomputer using  approximately $300$  Intel Xeon E5-2692 cores.
For the case with $512^3$ cells, the computational time of CPU code
with supercomputer using 1024 Intel Xeon E5-2692 cores is
approximately $3$ times longer than that of GPU code using $8$ Tesla
V100 GPUs. It can be inferred that the efficiency of GPU code with
$8$ Tesla V100 GPUs approximately equals to that of MPI code with
$3000$ supercomputer cores, which agree well with the multiple-GPU
accelerated finite difference code for multiphase flows
\cite{GPU-3}. These comparisons shows the excellent performance of
multiple-GPU accelerated HGKS for large-scale turbulence simulation.
To further show the performance of GPU computation, the scalability
is defined as
\begin{align}
S_{n}=\frac{T_{n,GPU}^{total}}{T_{1,GPU}^{total}}.
\end{align}
Fig.\ref{mpi-efficience} shows the log-log plot for $n$  and $S_n$,
where $2$, $4$ and $8$ GPUs are used for the case with $256^3$
cells, while $8$, $12$ and $16$ GPUs are used for the case with
$512^3$ cells. The ideal scalability of parallel computations would
be equal to $n$. With the log-log plot for $n$ and
$T_{n,GPU}^{total}$, an ideal scalability would follow $-1$ slope.
However, such idea scalability is not possible due to the
communication delay among the computational cores and the idle time
of computational nodes associated with load balancing. As expected,
the explicit formulation of HGKS scales properly with the increasing
number of GPU. Conceptually, the total computation amount increases
with a factor of $16$, when the number of cells doubles in every
direction. Taking the communication delay into account, the
execution time of  $256^3$ cells with 1 GPU approximately equals to
that of  $512^3$ cells with 16 GPUs, which also indicates the
scalability of GPU code. When GPU code using more than $8$ GPUs, the
communication across GPU nodes with RoCE is required, which accounts
for the worse scalability using $12$ GPUs and $16$ GPUs.
Table.\ref{GPU-CPU-TG-D} shows the execution time in flow solver
$T_{n,GPU}^{flow}$ and the communication delay between CPU and GPU
$T_{n,GPU}^{com}$.   Specifically, $T_{n,GPU}^{flow}$ is consist of
the time of WENO reconstruction, flux calculation at cell interface
and update of flow variables, while $T_{n,GPU}^{com}$ concludes time
for MPI$\_$RECV and MPI$\_$SEND for initialization and boundary
exchange, and MPI$\_$REDUCE for global time step. The histogram of
$T_{n,GPU}^{com}$ and $T_{n,GPU}^{flow}$ with different number of
GPU is shown in Fig.\ref{mpi-efficience-2} with $256^3$ and $512^3$
cells. With the increase of GPU, more time for communication is
consumed and the parallel efficiency is reduced accounting for the
practical scalability in Fig.\ref{mpi-efficience}. Especially, the
communication across GPU nodes with RoCE is needed for the GPU code
using more than $8$ GPUs, which consumes longer time than the
communication in single GPU node with NVLink. The performance of
communication across GPU nodes with InfiniBand will be tested, which
is designed for HPC center to achieve larger throughout among GPU
nodes.

\begin{figure}[!h]
\centering
\includegraphics[width=0.45\textwidth]{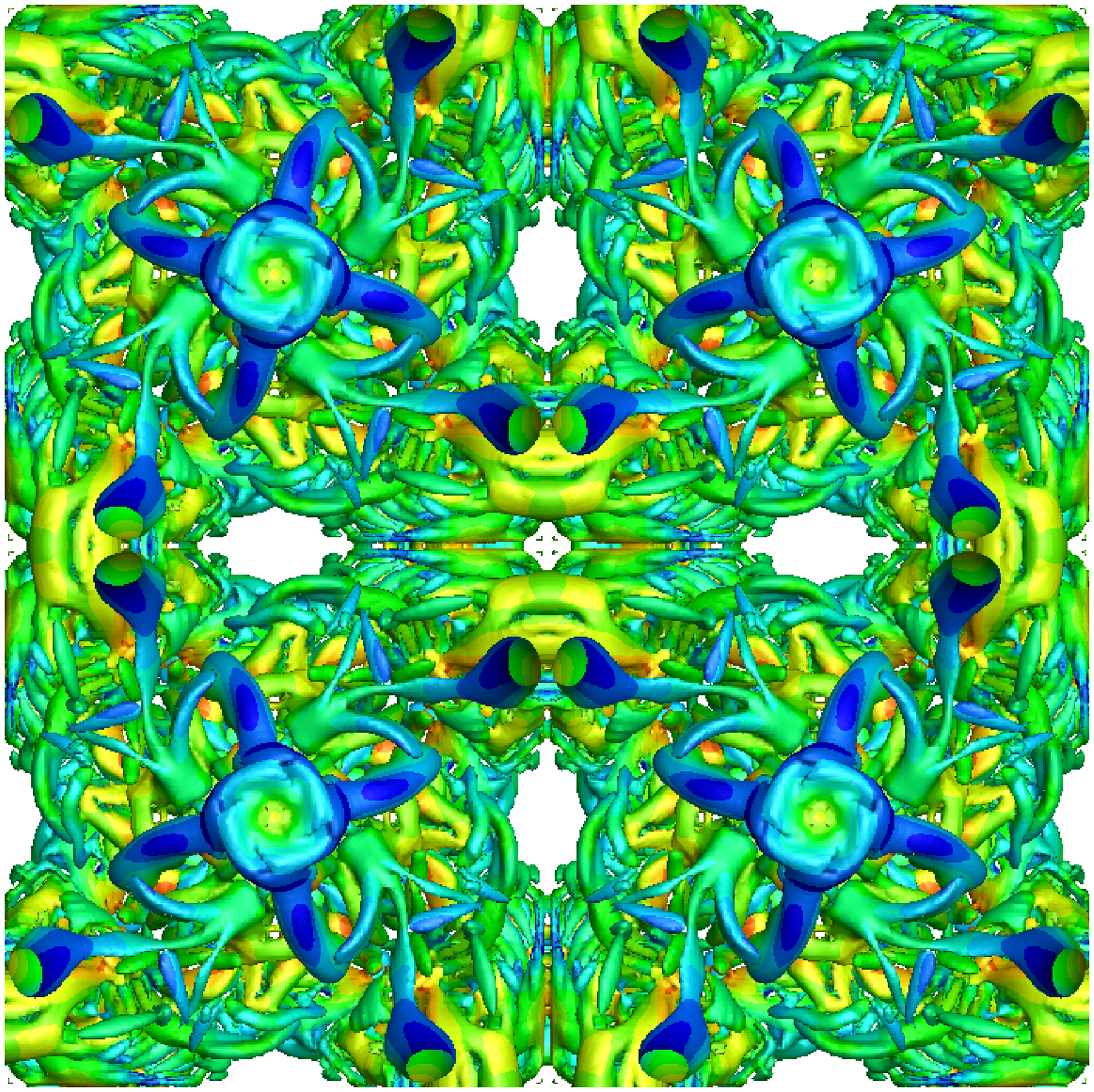}
\includegraphics[width=0.45\textwidth]{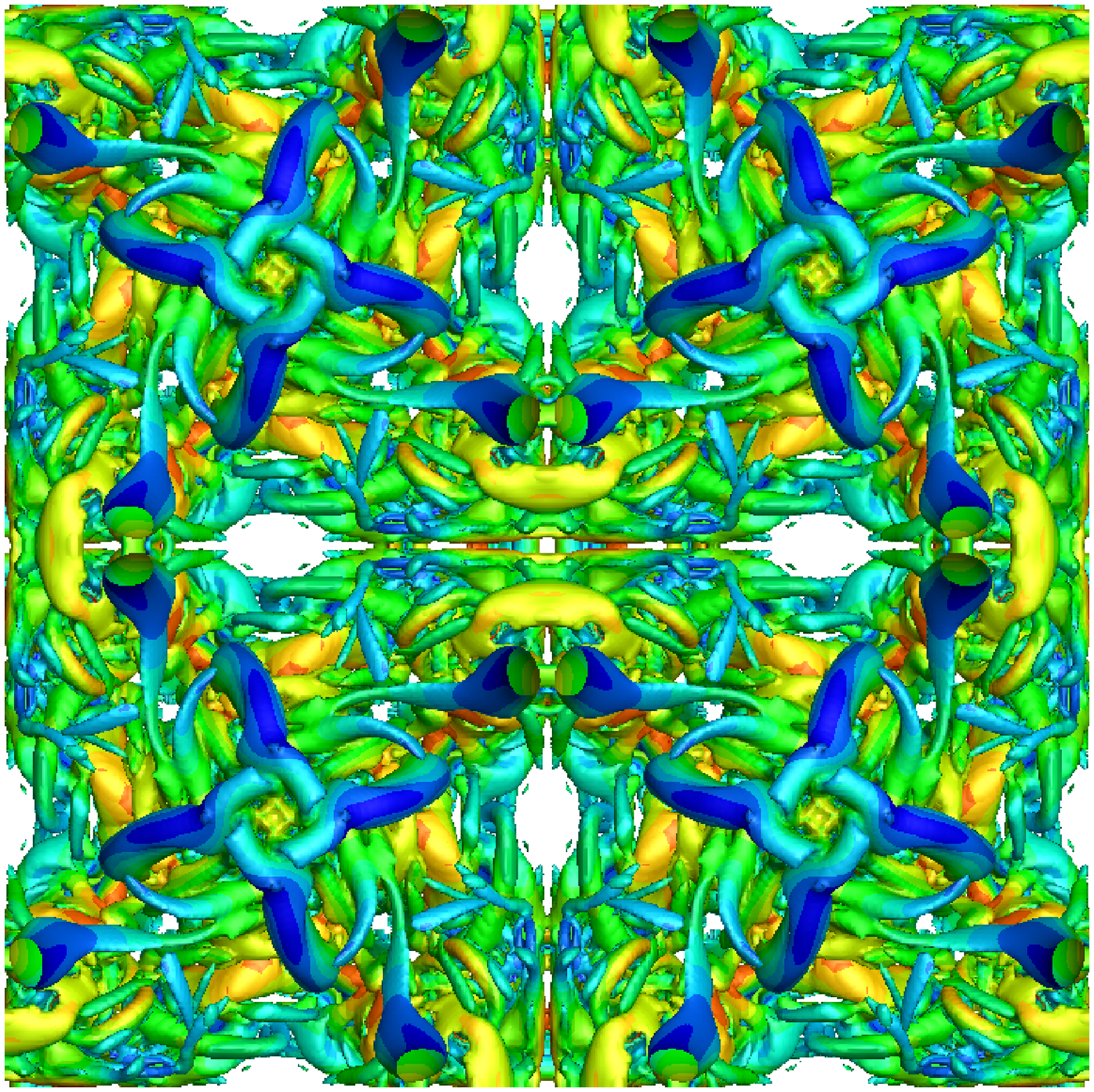}
\includegraphics[width=0.45\textwidth]{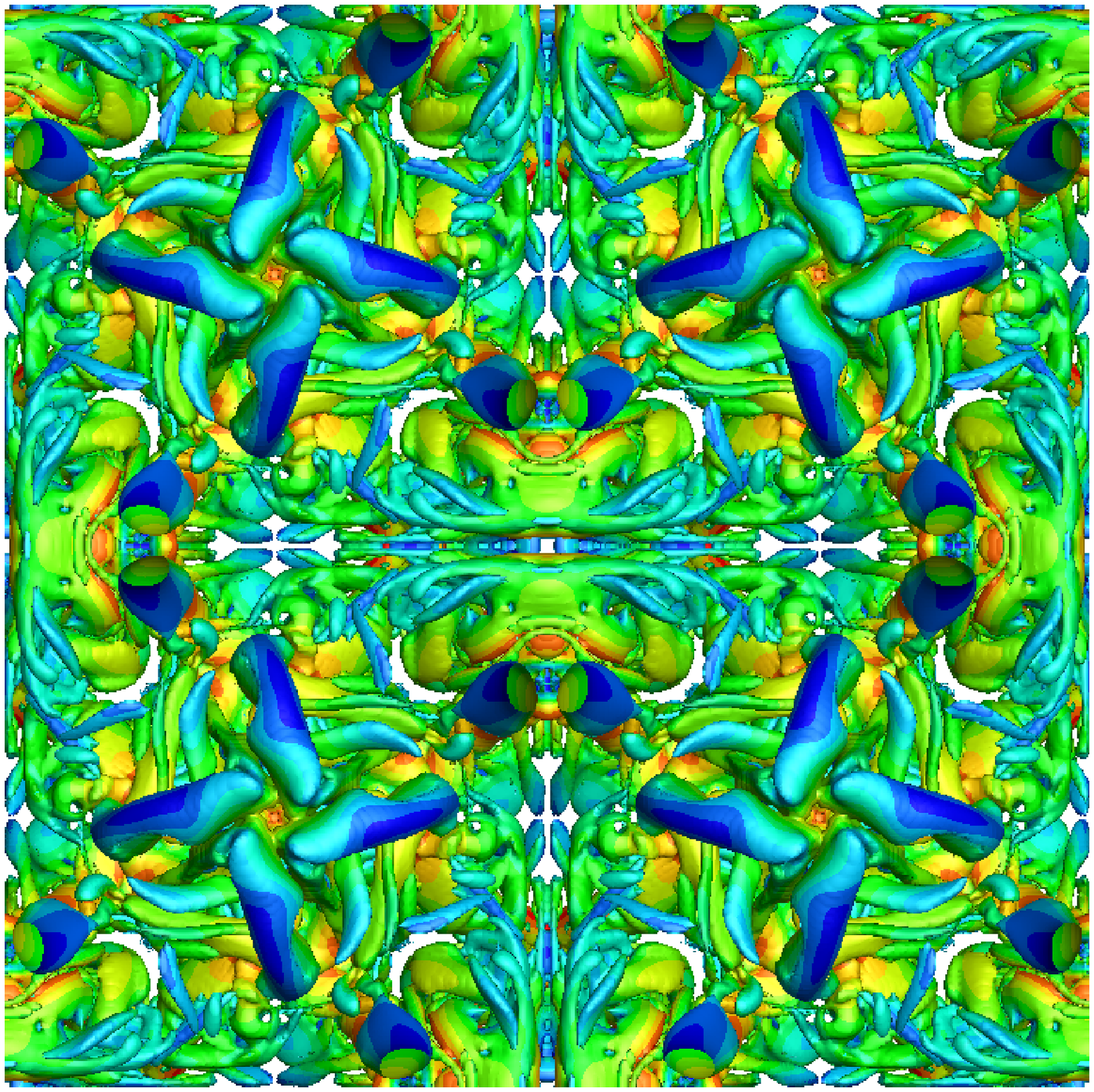}
\includegraphics[width=0.45\textwidth]{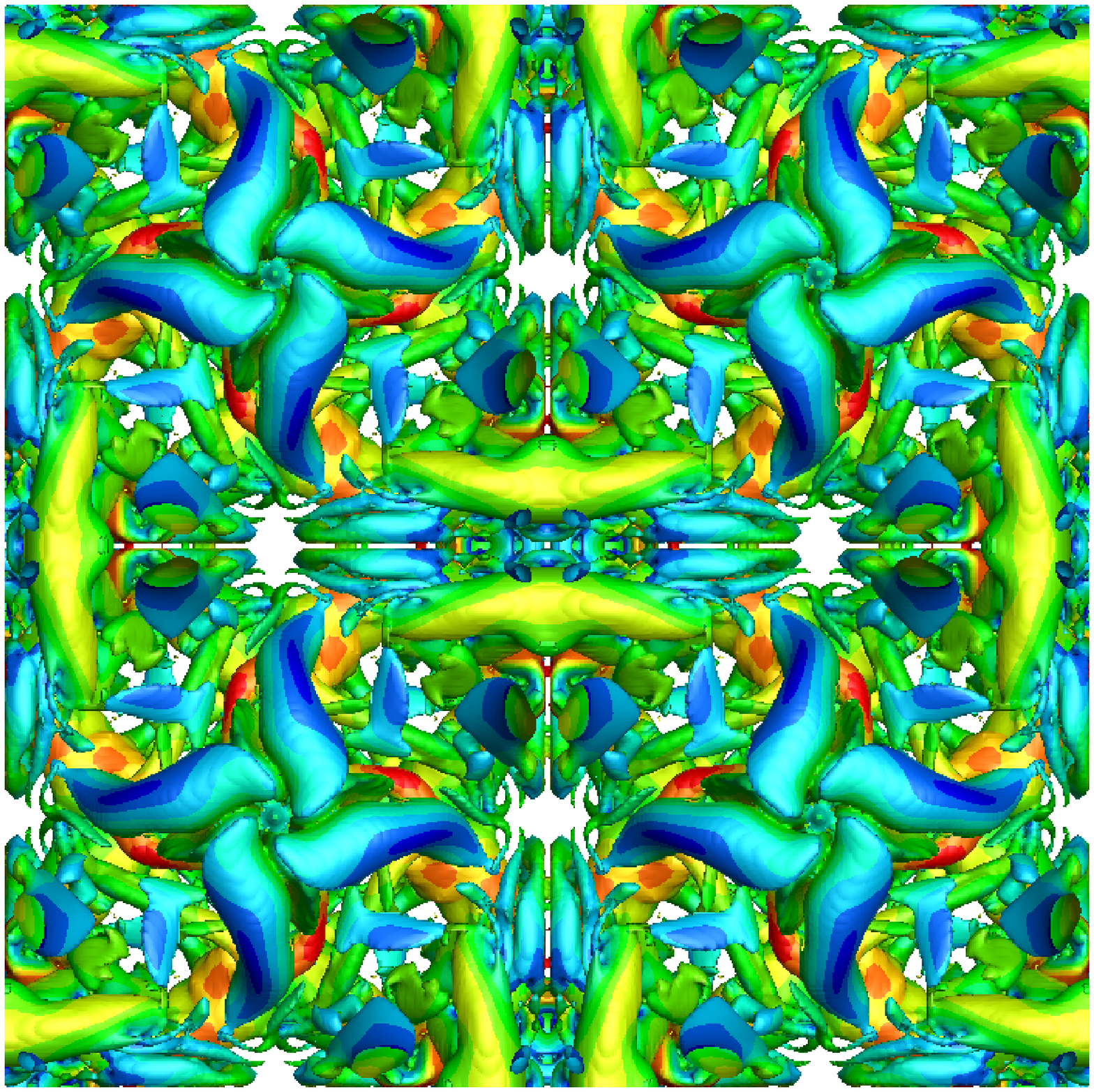}
\caption{\label{tg-vortex-q} Taylor-Green vortex: top view of
Q-criterion with Mach number $Ma = 0.25$, $0.5$, $0.75$ and $1.0$ at
$t =10$ using $256^3$ cells.}
\end{figure}

\begin{figure}[!h]
\centering
\includegraphics[width=0.495\textwidth]{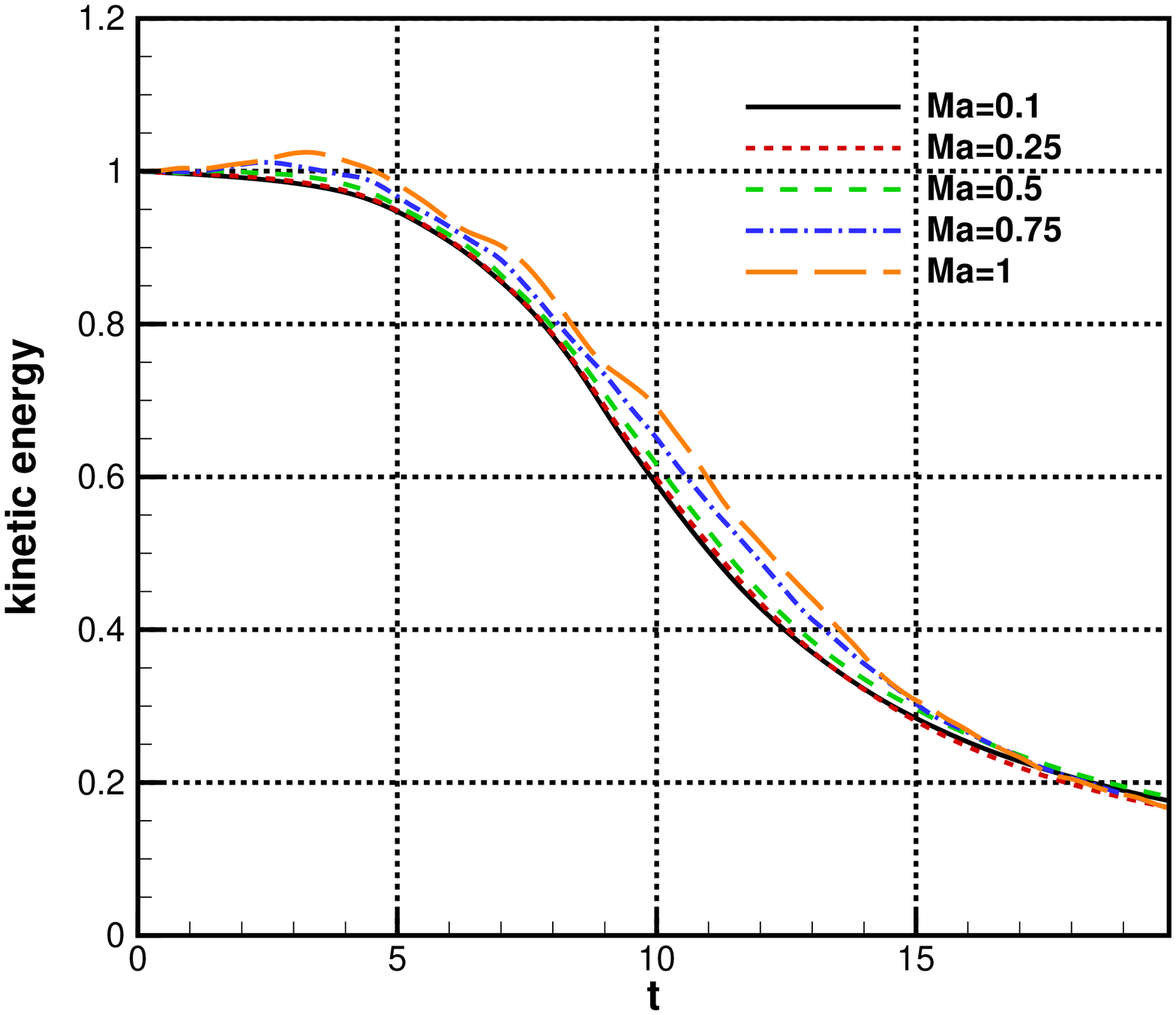}
\includegraphics[width=0.495\textwidth]{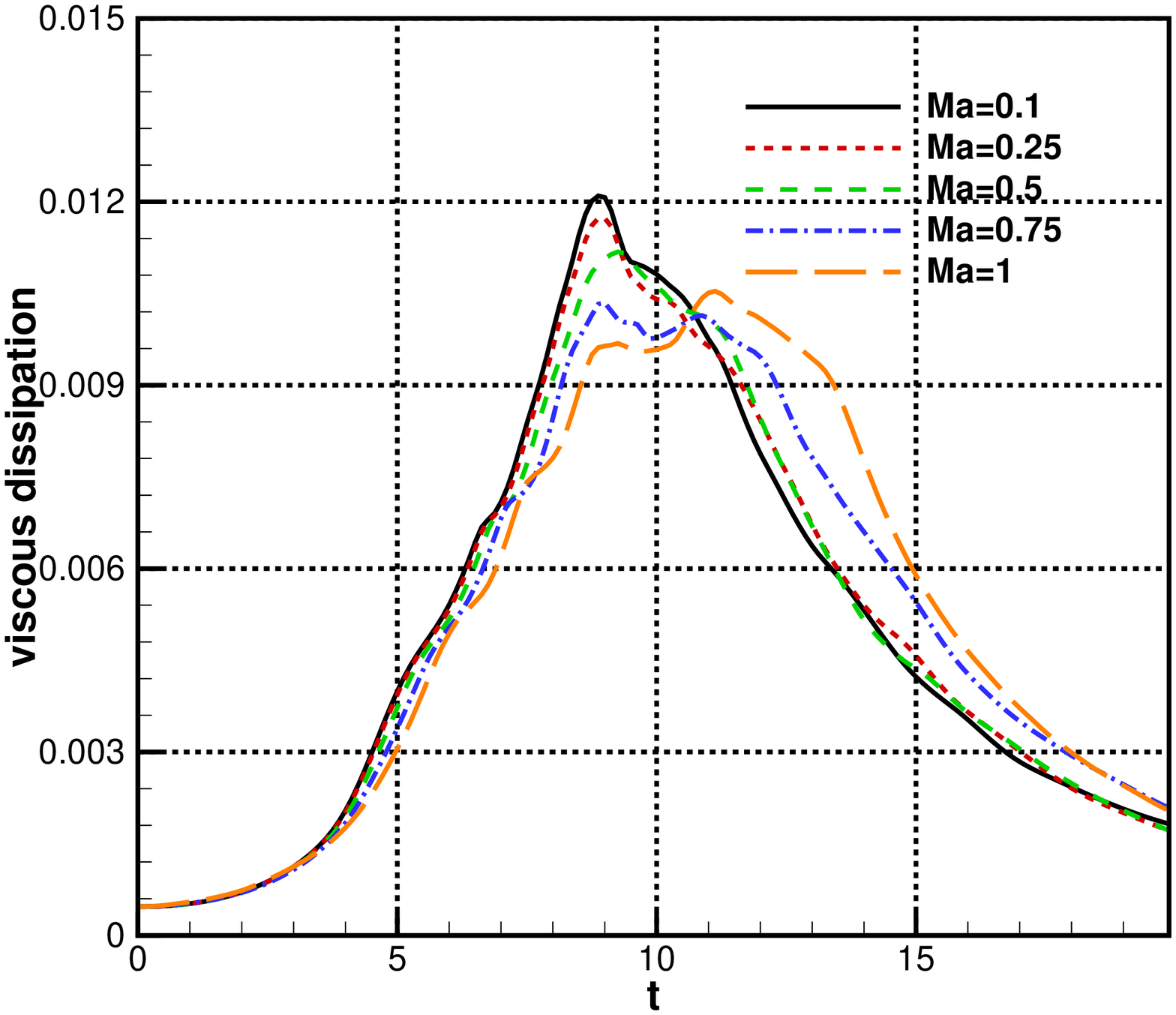}
\caption{\label{tg-vortex-compressible} Taylor-Green vortex: the
time history of kinetic energy  $E_k$ and viscous dissipation
$\varepsilon_{com}$ with Mach number $0.1$, $0.25$, $0.5$, $0.75$
and $1.0$ using $256^3$ cells.}
\end{figure}

The time history of kinetic energy $E_k$, dissipation rate
$\varepsilon(E_k)$ and enstrophy dissipation rate
$\varepsilon(\zeta)$ of the above nearly incompressible case from
GPU code is identical with the solution from previous CPU code, and
more details can be found in \cite{GKS-high-DNS}. With the efficient
multiple-GPU accelerated platform, the compressible Taylor-Green
vortexes with fixed Reynolds number $Re = 1600$ Mach number $Ma =
0.25$, $0.5$, $0.75$ and $1.0$ are also tested. Top view of
Q-criterion with Mach number $Ma = 0.25$, $0.5$, $0.75$ and $1.0$ at
$t =10$ using $256^3$ cells are presented in Fig.\ref{tg-vortex-q}.
The volume-averaged kinetic energy is given by
\begin{align*}
E_k=\frac{1}{\rho_0\Omega}\int_\Omega\frac{1}{2}\rho\boldsymbol{U}\cdot\boldsymbol{U} \text{d} \Omega,
\end{align*}
where $\Omega$ is the volume of the computational domain.
The total viscous dissipation is defined as
\begin{align*}
    \varepsilon_{com}=\frac{\mu}{\rho_0\Omega}\int_\Omega\boldsymbol{\omega}\cdot\boldsymbol{\omega} \text{d} \Omega
    +\frac{4}{3}\frac{\mu}{\rho_0\Omega}\int_\Omega(\nabla\cdot\boldsymbol{U})^2 \text{d} \Omega,
\end{align*}
where $\mu$ is the coefficient of viscosity and
$\boldsymbol{\omega}=\nabla\times \boldsymbol{U}$. The time
histories of kinetic energy and total viscous dissipation are given
in Fig.\ref{tg-vortex-compressible}. Compared with the low Mach
number cases, the transonic and supersonic cases with $Ma = 0.75$
and $Ma = 1.0$ have regions of increasing kinetic energy rather than
a monotone decline. This is consistent with the compressible TGV
simulations \cite{Case-Peng}, where increases in the kinetic energy
were reported during the time interval  $2 \leq t \leq 4$. A large
spread of total viscous dissipation curves is observed in
Fig.\ref{tg-vortex-compressible}, with the main trends being a delay
and a flattening of the peak dissipation rate for increasing Mach
numbers. The peak viscus dissipation rate for the $Ma=0.25$ curve
matches closely to that of $Ma=0.1$. Meanwhile, for the higher Mach
number cases, the viscous dissipation rate is significantly higher
in the late evolution region (i.e., $t \ge 11$) than that from near
incompressible case.

\subsection{Turbulent channel flows for FP32-FP64 comparison}
\begin{figure}[!h]
\centering
\includegraphics[width=0.475\textwidth]{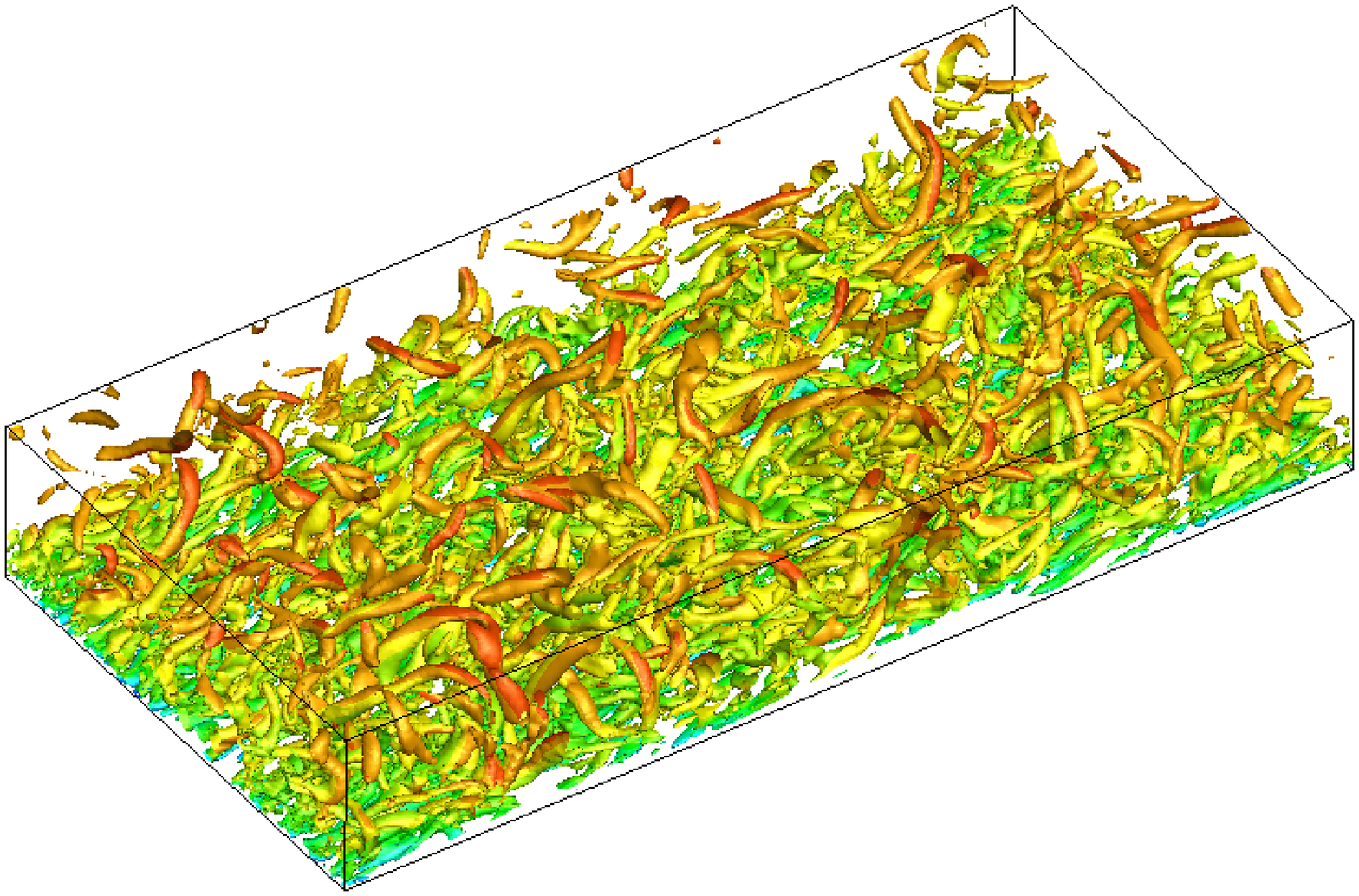}
\includegraphics[width=0.475\textwidth]{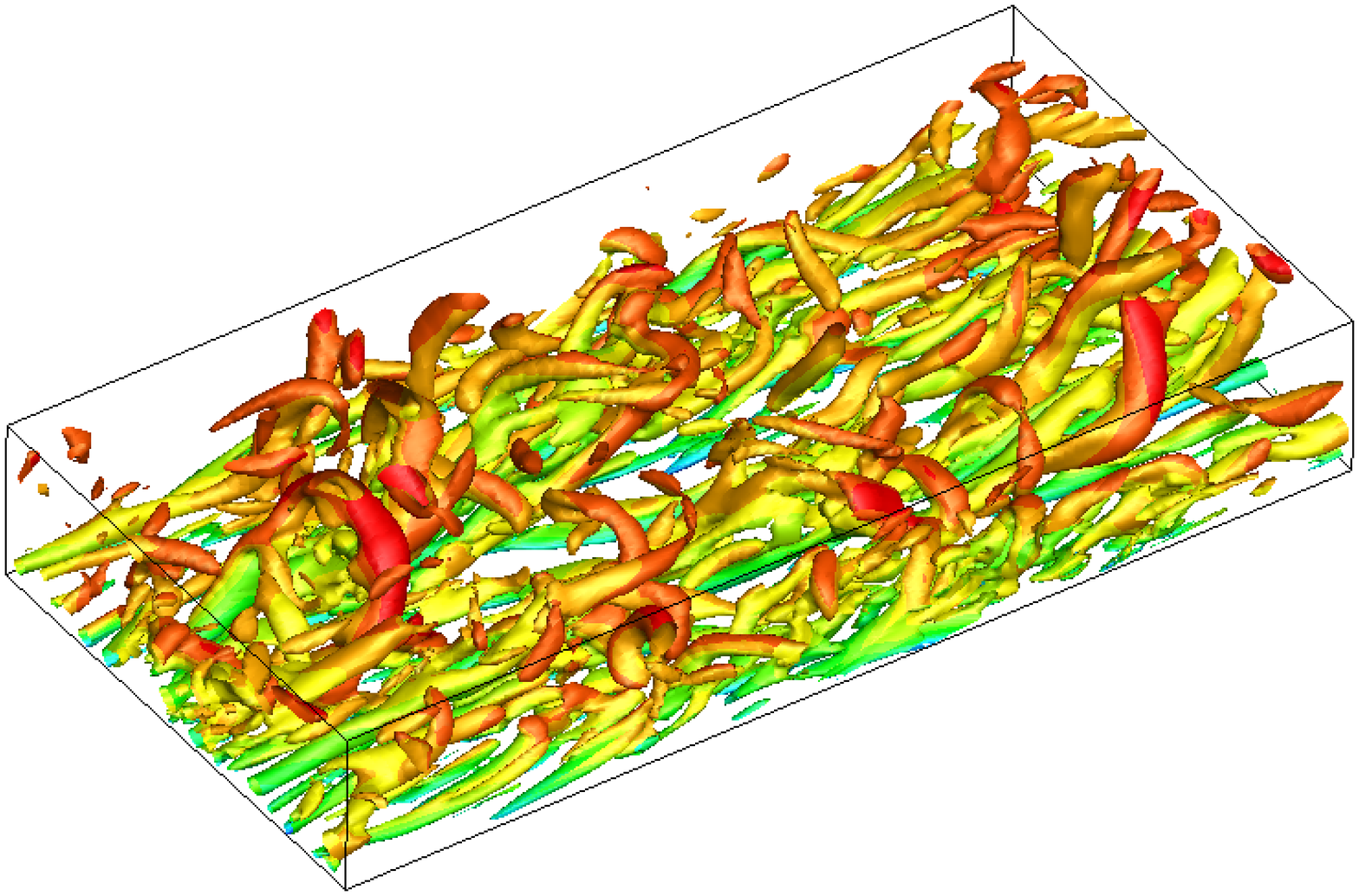}
\caption{\label{channel_iso_surface} Turbulent channel flow: the
instantaneous $Q$-criterion iso-surfaces with $G_2$ (left) and $H_2$
(right).}
\end{figure}

For most GPUs, the  performance of FP32  precision is stronger than
FP64 precision. In this case, to address the effect of FP32 and FP64
precision in terms of efficiency, memory costs and simulation
accuracy, both the near incompressible and compressible turbulent
channel flows are tested. Incompressible and compressible turbulent
channel flows have been studied to understand the mechanism of
wall-bounded turbulent flows
\cite{kim1987turbulence,kim1999turbulence,coleman1995numerical}. In
the computation, the physical domain is $(x,y,z)\in[0,2 \pi H]
\times [-H,H] \times[0,\pi H]$, and the computational mesh is given
by the coordinate transformation
\begin{align*}
\begin{cases}
\displaystyle x=\xi,\\
\displaystyle y=\tanh(b_g(\frac{\eta}{1.5\pi}-1))/\tanh(b_g),\\
\displaystyle z=\zeta,
\end{cases}
\end{align*}
where the computational domain takes $(\xi,\eta,\zeta)\in[0, 2\pi
H]\times[0, 3\pi H]\times[0, \pi H]$ and $b_g=2$. The fluid is
initiated with density $\rho = 1$ and the initial streamwise
velocity $U(y)$ profile is given by the perturbed Poiseuille flow
profile, where the white noise is added with $10\%$ amplitude of
local streamwise velocity. The spanwise and wall-normal velocity is
initiated with white noise. The periodic boundary conditions are
used in streamwise $x$-direction and spanwise $z$-directions, and
the non-slip and isothermal boundary conditions are used in vertical
$y$-direction. The constant moment flux in the streamwise direction
is used to determine the external force \cite{GKS-high-DNS}. In
current study, the nearly incompressible turbulent channel flow with
friction Reynolds number $Re_\tau=395$ and Mach number $Ma=0.1$, and
the compressible turbulent channel flow with bulk Reynolds number
$Re=3000$ and bulk Mach number $Ma=0.8$ are simulated with the same
set-up as previous studies \cite{GKS-high-DNS,GKS-high-cao-iles}.
For compressible case, the viscosity $\mu$ is determined by the
power law as $\mu= \mu_w(T/T_w)^{0.7}$, and the Prandtl number $Pr
=0.7$ and $\gamma=1.4$ are used. For the nearly incompressible
turbulent channel flow, two cases $G_1$ and $G_2$ are tested, where
$256\times128\times128$ and $256\times256\times128$ cells are
distributed uniformly in the computational space as
Table.\ref{channel_grids-1}. For the compressible flow, two cases
$H_1$ and $H_2$ are tested as well, where $128\times128\times128$
and $128\times256\times128$ cells are distributed uniformly in the
computational space as Table.\ref{channel_grids-2}. In the
computation, the cases $G_1$ and $H_1$ are implemented with single
GPU, while the cases $G_2$ and $H_2$ are implemented with double
GPUs. Specifically, $\Delta y^{+}_{min}$ and $\Delta y^{+}_{max}$
are  the minimum and maximum mesh space in the $y$-direction. The
equidistant meshes are used in $x$ and $z$ directions, and $\Delta
x^{+}$ and $\Delta z^{+}$ are the equivalent plus unit,
respectively. The fine meshes arrangement in case $G_2$ and $H_2$
can meet the mesh resolution for DNS in turbulent channel flows. The
$Q$-criterion iso-surfaces for $G_2$ and $H_2$ are given
Fig.\ref{channel_iso_surface} in which both cases can well resolve
the vortex structures. With the mesh refinement, the high-accuracy
scheme corresponds to the resolution of abundant turbulent
structures.

\begin{table}[!htp]
\centering
\begin{tabular}{c|c|c|c|c}
\hline
\hline
Case    &Mesh size  &$\Delta y^{+}_{min}$/$\Delta y^{+}_{max}$   &$\Delta x^{+}$   &$\Delta z^{+} $\\
\hline
$G_1$   &$256\times128\times128$   &0.93/12.80   &9.69   &9.69 \\
\hline
$G_2$   &$256\times256\times128$  &0.45/6.40   &9.69    &9.69 \\
\hline
\hline
\end{tabular}
\caption{\label{channel_grids-1} Turbulent channel flow: the details
of mesh for nearly incompressible flow with friction Reynolds number
$Re_\tau=395$ and $Ma=0.1$.}
~\\
\centering
\begin{tabular}{c|c|c|c|c}
\hline
\hline
Case    &Mesh size  &$\Delta y^{+}_{min}$/$\Delta y^{+}_{max}$   &$\Delta x^{+}$   &$\Delta z^{+} $\\
\hline
$H_1$   &$128\times128\times128$   &0.43/5.80   &8.76   &4.38 \\
\hline
$H_2$   &$128\times256\times128$  &0.22/2.90   &8.76   &4.38 \\
\hline
\hline
\end{tabular}
\caption{\label{channel_grids-2} Turbulent channel flow: the
details of mesh for compressible flow with bulk Reynolds number
$Re=3000$ and bulk Mach number $Ma=0.8$.}
\end{table}

\begin{table}[!h]
\centering
\begin{tabular}{c|c|c|c}
\hline
\hline
TITAN  RTX   & $T_{2,GPU}^{total}$ with FP32  &  $T_{2,GPU}^{total}$ with FP64  & $S_{fp}$   \\
\hline
$G_1$     &0.414    &1.152  & 2.783   \\
\hline
$G_2$    &1.635      &4.366   & 2.670  \\
\hline
TITAN  RTX   &  memory cost with FP32  &  memory cost with FP64 & $R_{fp}$  \\
\hline
$G_1$     &2.761    &5.355   & 1.940   \\
\hline
$G_2$    &5.087     &10.009   & 1.968  \\
\hline
\hline
Tesla V100    & $T_{2,GPU}^{total}$ with FP32 &  $T_{2,GPU}^{total}$ with FP64  & $S_{fp}$   \\
\hline
$G_1$     &0.377    & 0.642  & 1.703   \\
\hline
$G_2$    &1.484      &2.607   & 1.757  \\
\hline
Tesla V100   &  memory cost with FP32   &  memory cost with FP64  & $R_{fp}$  \\
\hline
$G_1$     &3.064    &5.906   & 1.928   \\
\hline
$G_2$    &5.390      &10.559   & 1.959  \\
\hline
\hline
\end{tabular}
\caption{\label{channel-precision} Turbulent channel flow: the
comparison of computational times and memory for FP32 and FP64
precision.}
\end{table}

Because of the reduction in device memory and improvement of
arithmetic capabilities on GPUs, the benefits can be achieved by
using FP32 precision compared to FP64 precision. In view of these
strength, FP32-based and mixed-precision-based high-performance
computing start to be explored \cite{lehmann2021accuracy,
haidar2020mixed}. However, whether FP32 precision is sufficient for
the DNS of turbulent flows still needs to be discussed. It is
necessary to evaluate the effect of precision on DNS of compressible
turbulence. Thus, the GPU-accelerated HGKS is compiled with FP32
precision and FP64 precision, and both nearly incompressible and
compressible turbulent channel flows are used to test the
performance in different precision. For simplicity, only the memory
cost and computational time for the nearly incompressible cases are
provided in Table.\ref{channel-precision}, where the execution times
$T_{2,GPU}^{total}$ are given in terms of hours and the memory cost
is in GB. In the computation, two GPUs are used and $10H/U_c$
statistical period is simulated. In Table.\ref{channel-precision},
$S_{fp}$ represents the speed up of $T_{2,GPU}^{total}$ with FP32 to
$T_{2,GPU}^{total}$ with FP64, and $R_{fp}$ denotes the ratio of
memory cost with FP32 to memory cost with FP64. As expected, the
memory of  FP32 precision is about half of that of FP64 precision
for both TITAN RTX and Tesla V100 GPUs. Compared with FP64-based
simulation, 2.7x speedup is achieved for TITAN RTX GPU and 1.7x
speedup is achieved for Tesla V100 GPU.  As shown in
Table.\ref{channel-precision}, the comparable FP32 precision
performance between TITAN RTX and Tesla V100 also provides
comparable $T_{2,GPU}^{total}$ for these two GPUs. We conclude that
Turing architecture and Ampere architecture perform comparably for
memory-intensive computing tasks in FP32 precision.

\begin{figure}[!h]
\centering
\includegraphics[width=0.495\textwidth]{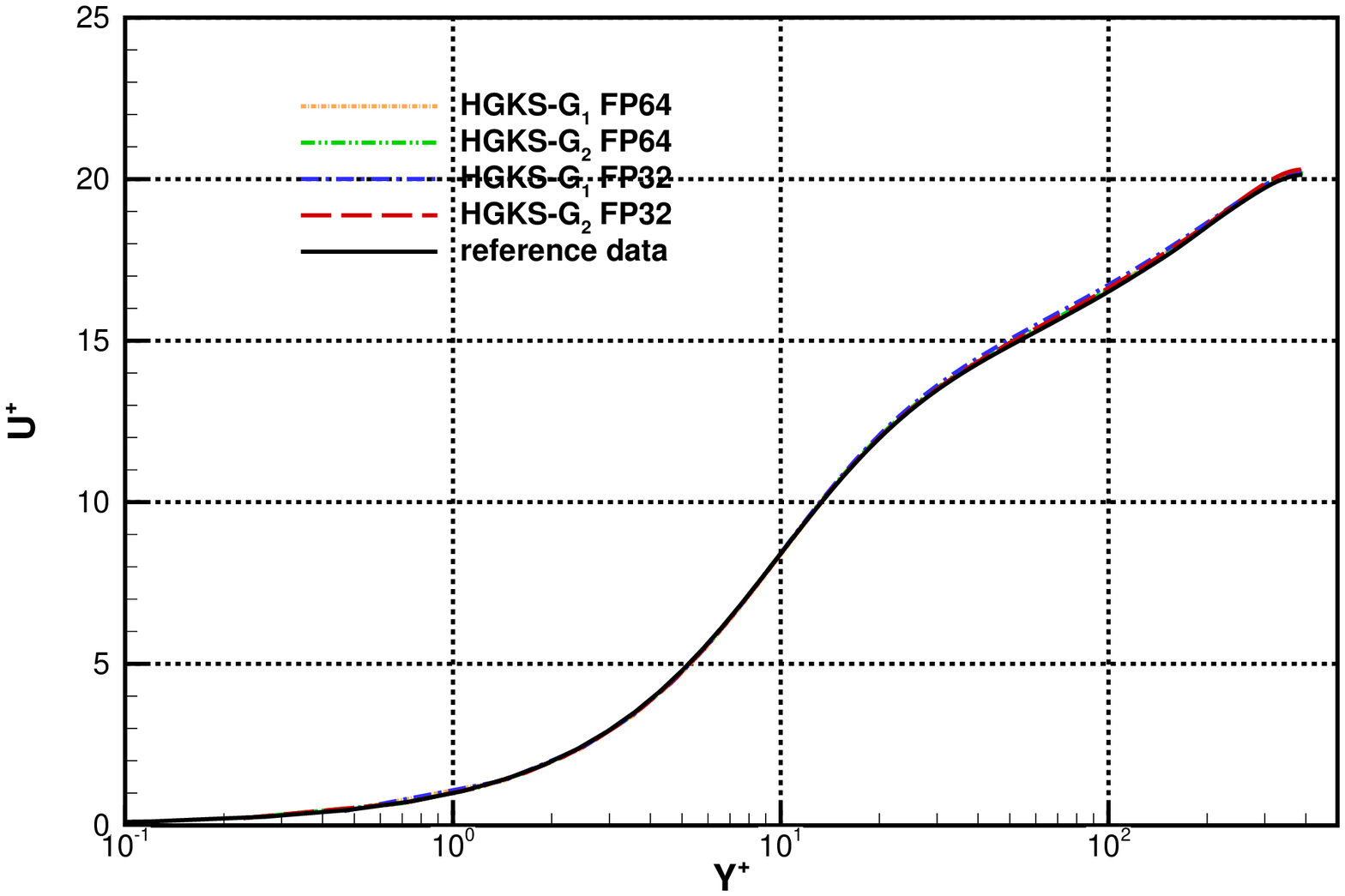}
\caption{\label{channel_ave_1} Turbulent channel flow: the mean
velocity for nearly incompressible turbulent flow. The reference
data is given in \cite{kim1999turbulence}.} \centering
\includegraphics[width=0.495\textwidth]{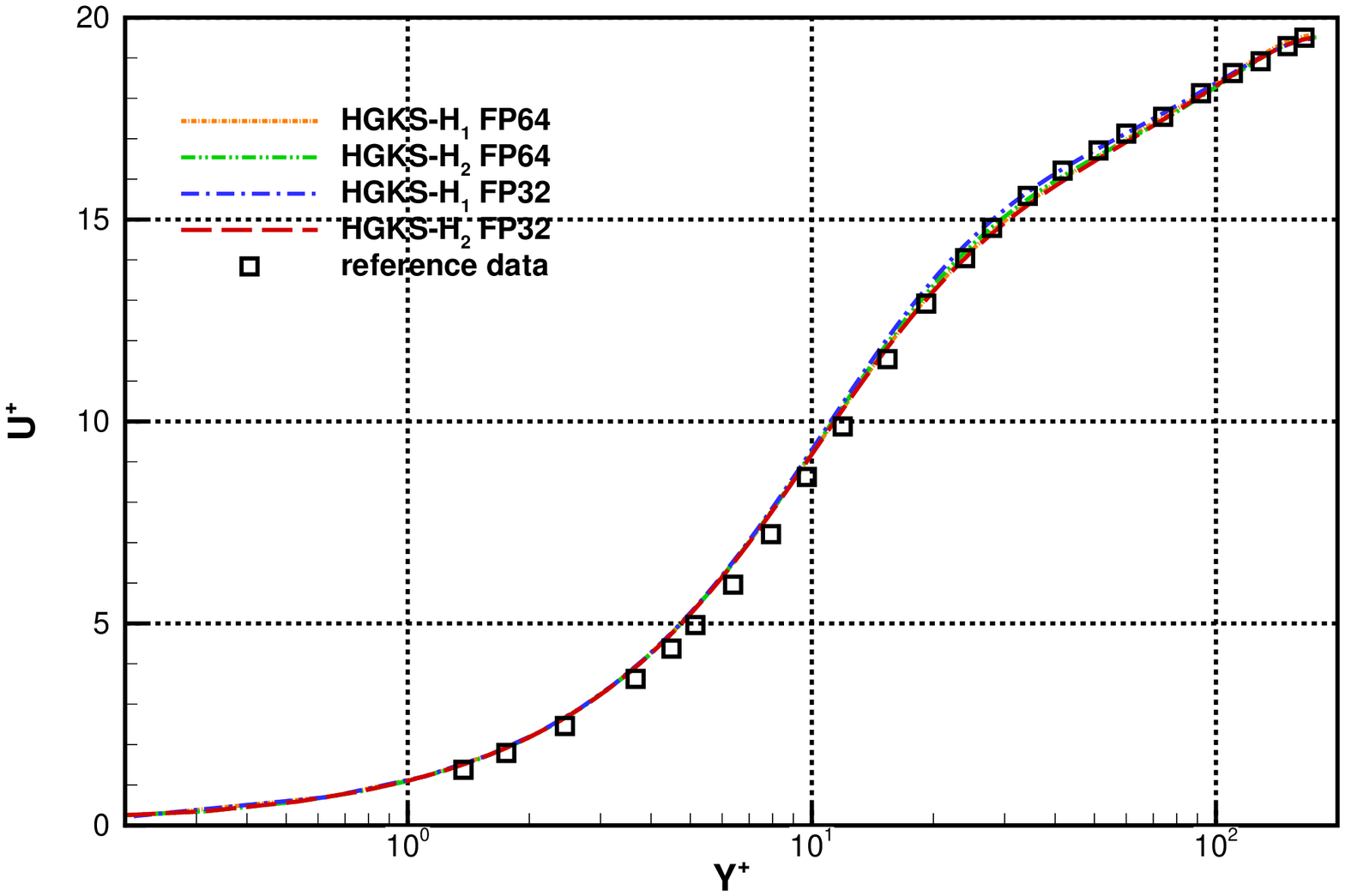}
\includegraphics[width=0.495\textwidth]{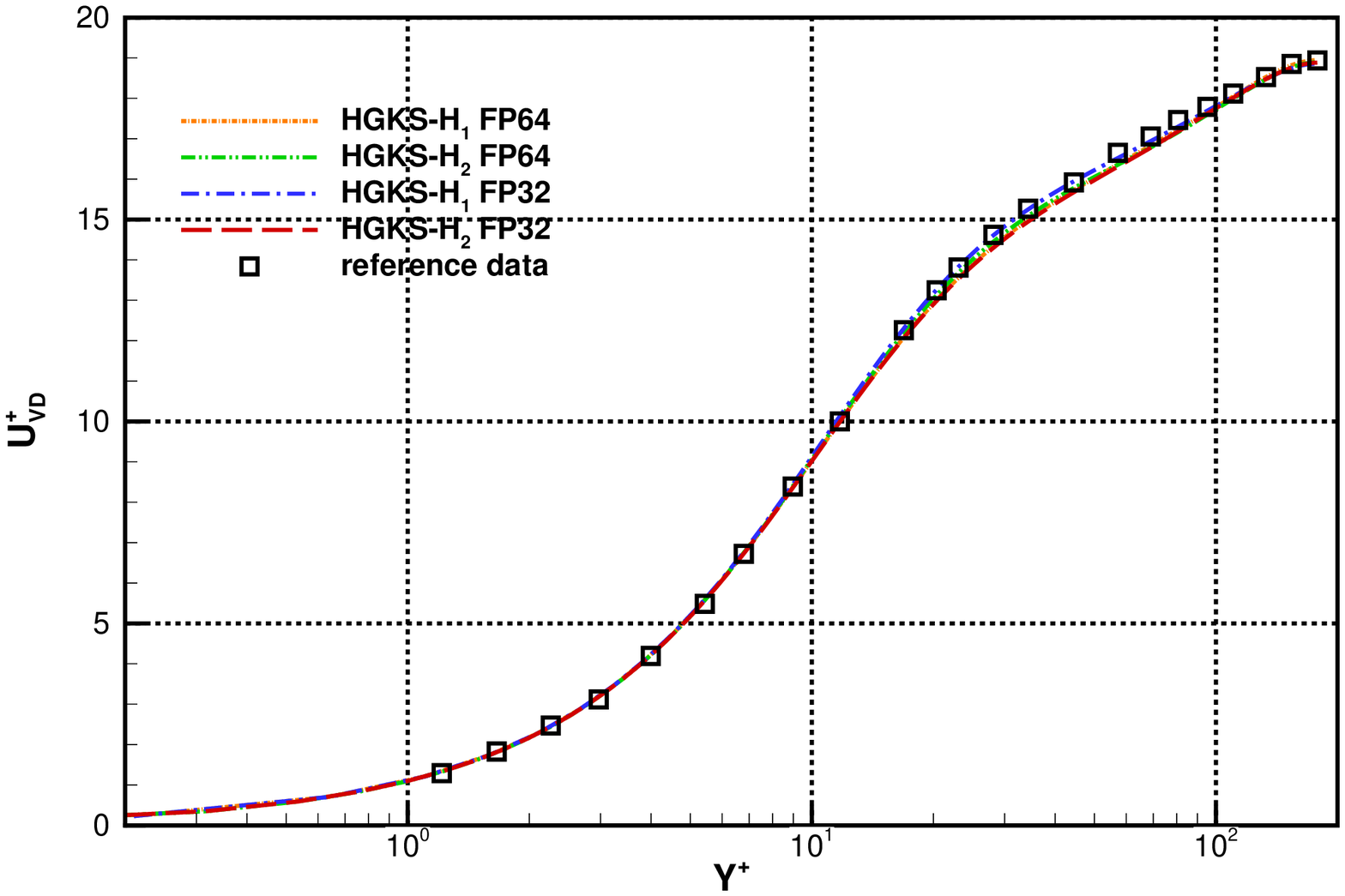}
\caption{\label{channel_ave_2} Turbulent channel flow: the mean
velocity $\langle U \rangle^+$ and velocity with VD transformation
$\langle U \rangle_{VD}^+$ for compressible turbulent flow. The
reference data is given in \cite{DNS-Li}.}
\end{figure}

\begin{figure}[!h]
\centering
\includegraphics[width=0.45\textwidth]{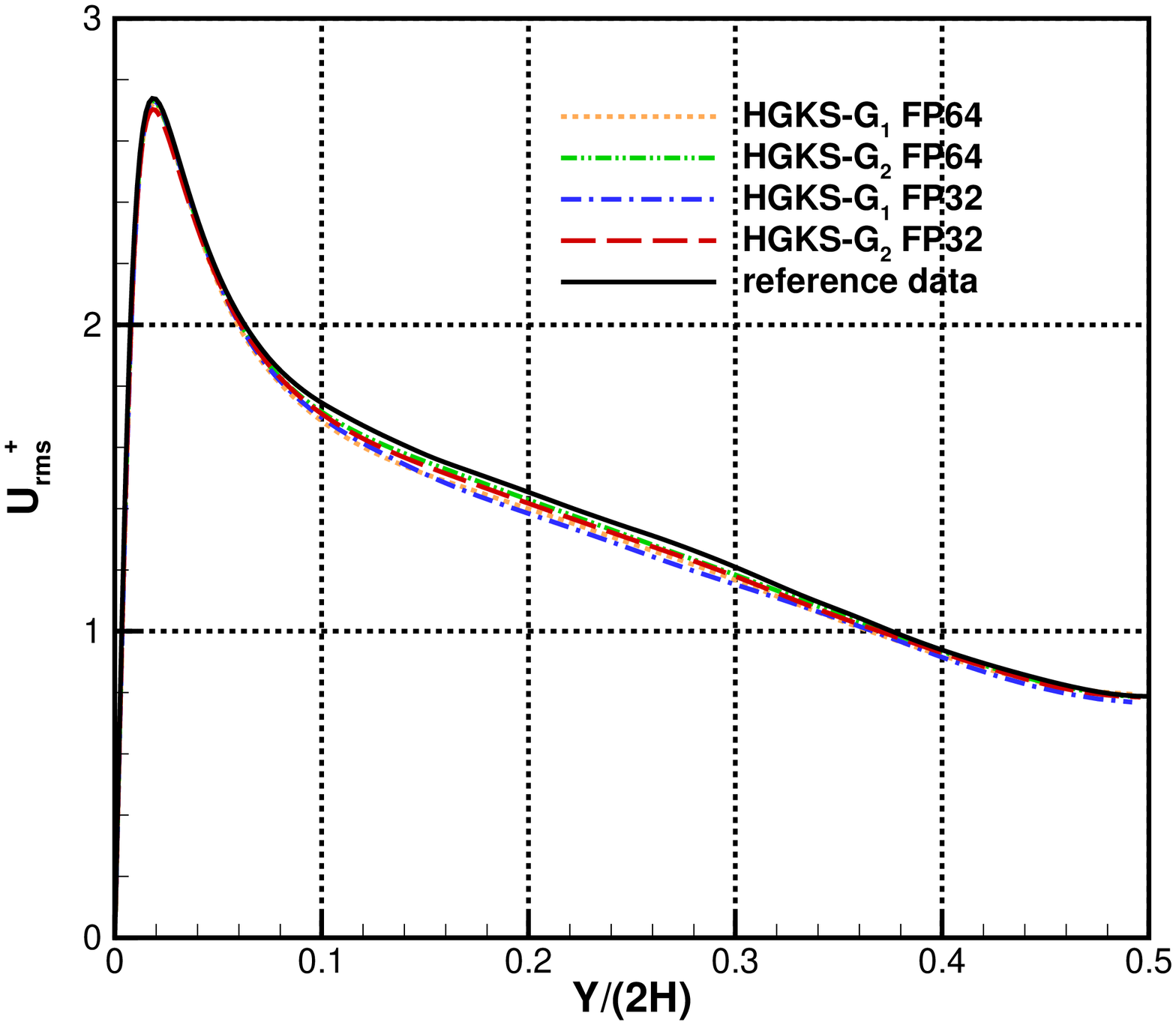}
\includegraphics[width=0.45\textwidth]{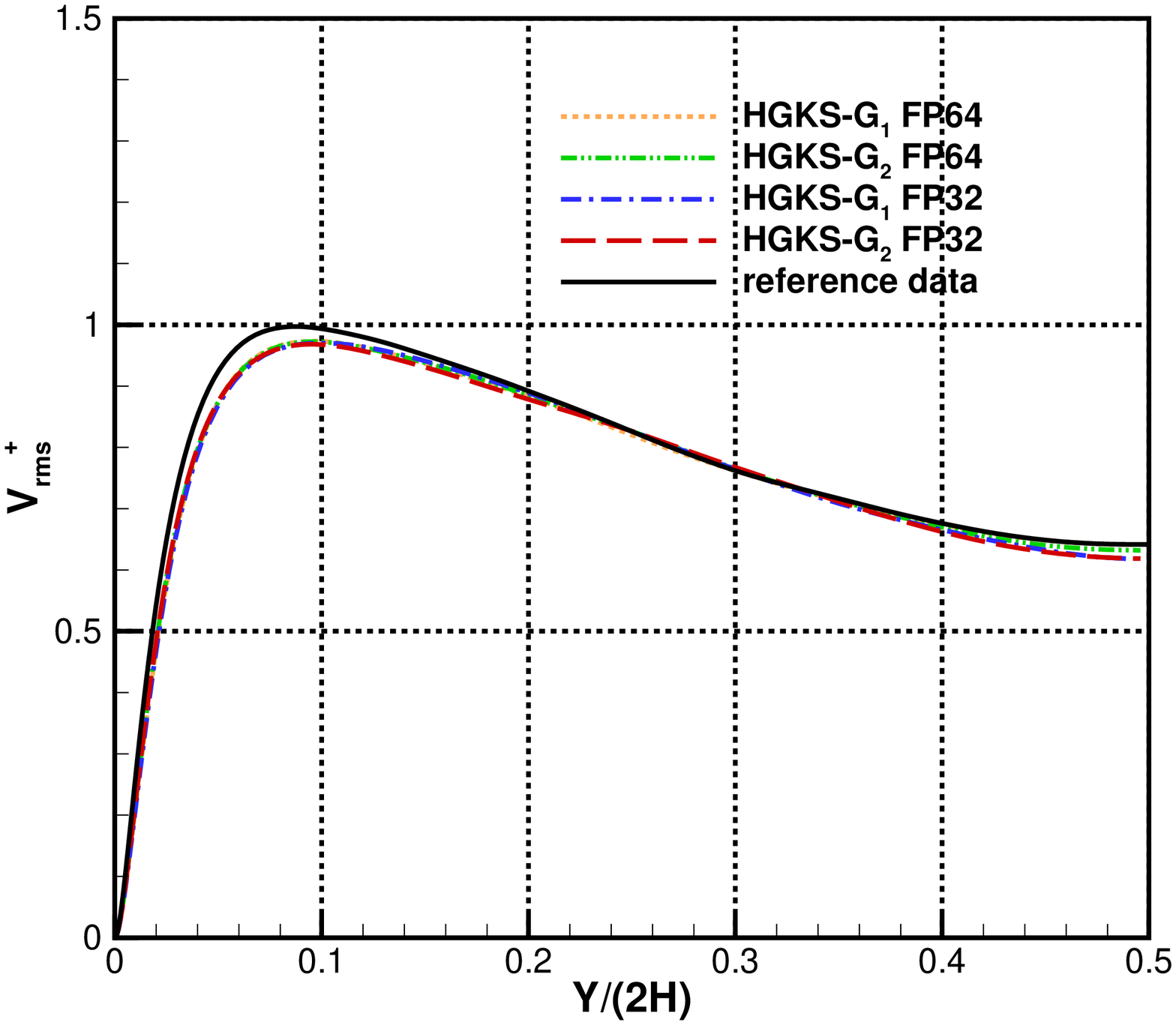}
\includegraphics[width=0.45\textwidth]{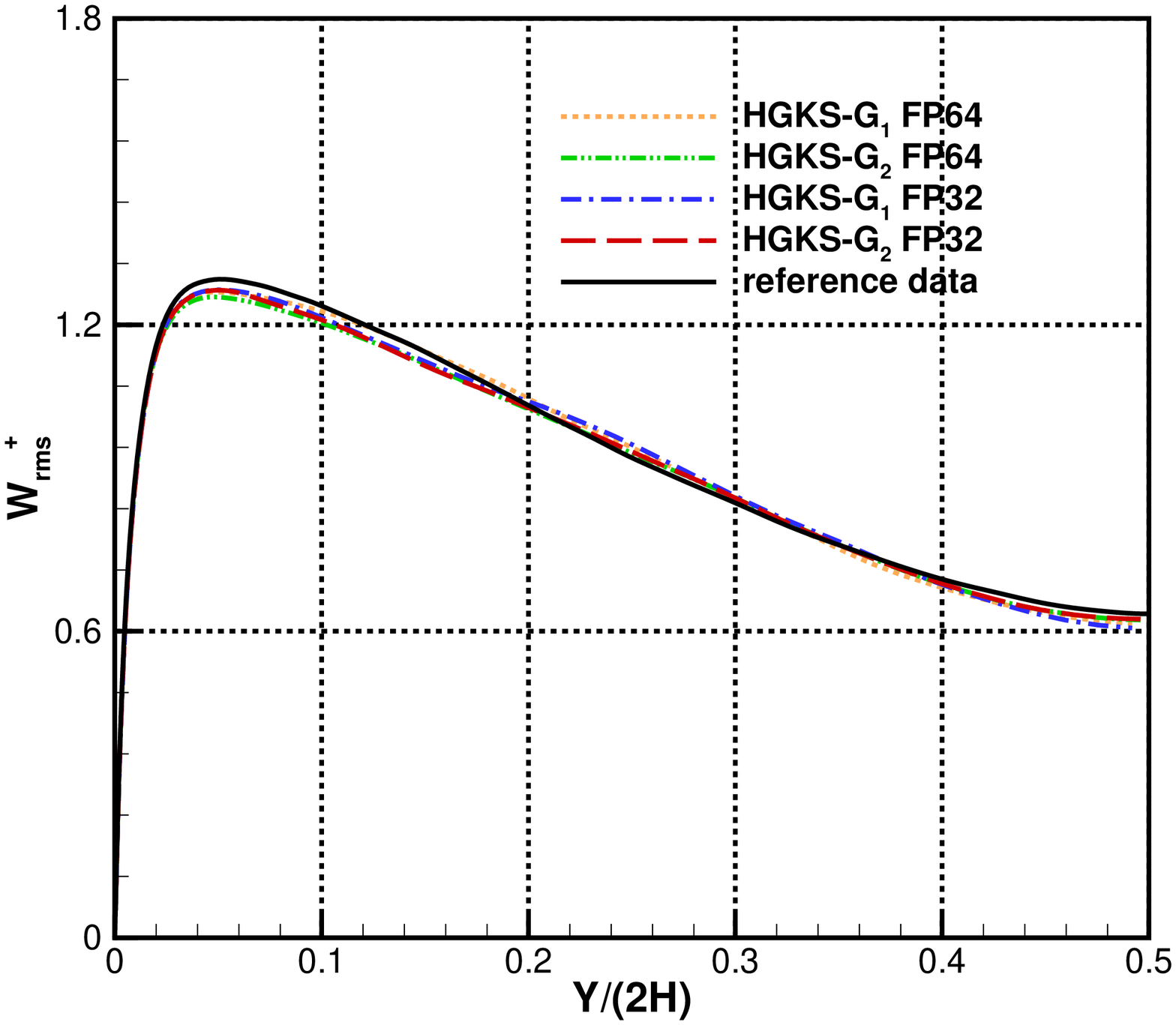}
\includegraphics[width=0.45\textwidth]{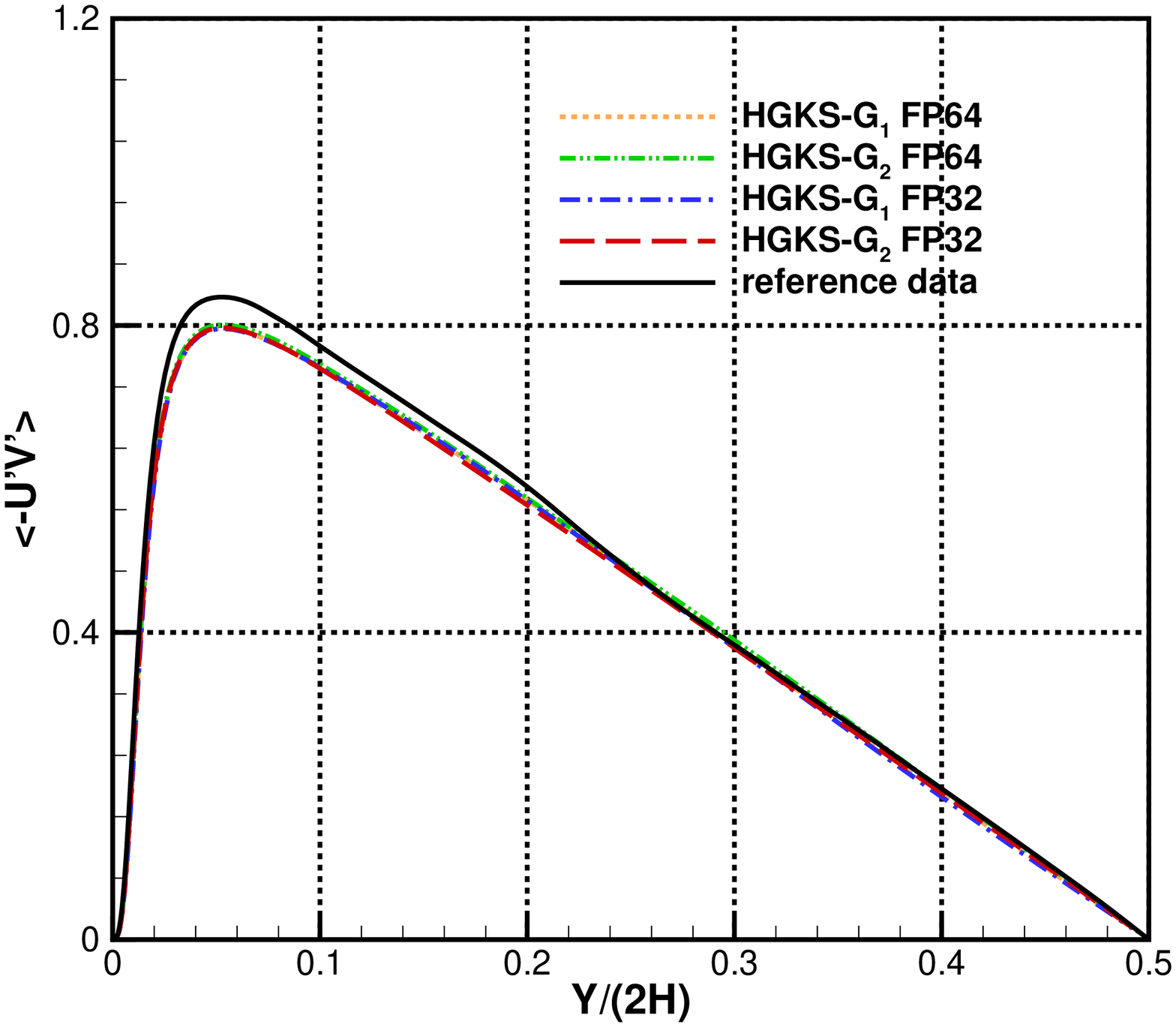}
\caption{\label{channel_fluctuation_1} Turbulent channel flow: the
root-mean-square fluctuation velocity $U_{rms}^+$, $V_{rms}^+$,
$W_{rms}^+$ and Reynolds stress $-<U'V'>$ profiles for nearly
incompressible turbulent flow. The reference data is given in
\cite{kim1999turbulence}.}
\end{figure}

\begin{figure}[!htp]
\centering
\includegraphics[width=0.45\textwidth]{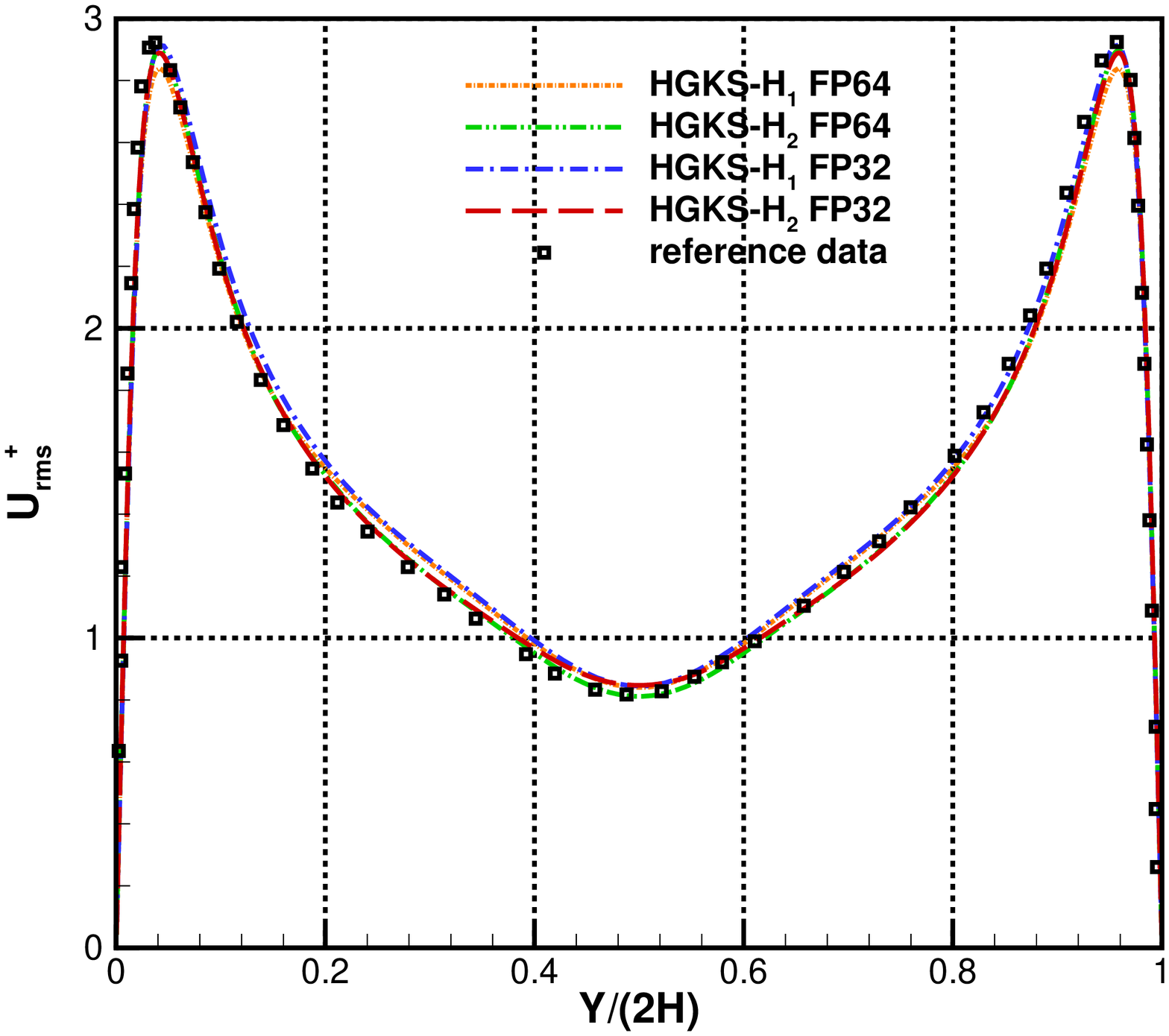}
\includegraphics[width=0.45\textwidth]{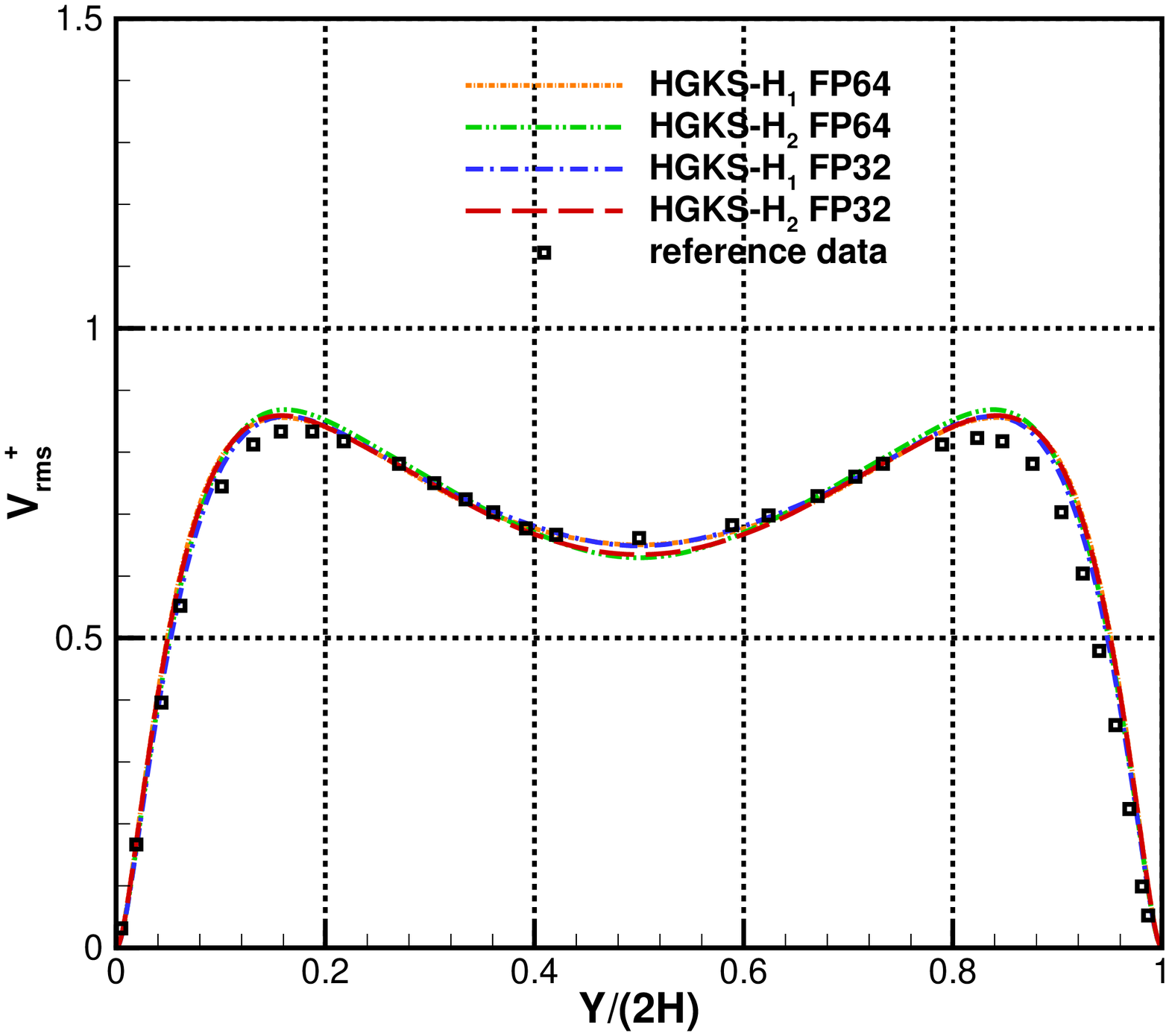}
\includegraphics[width=0.45\textwidth]{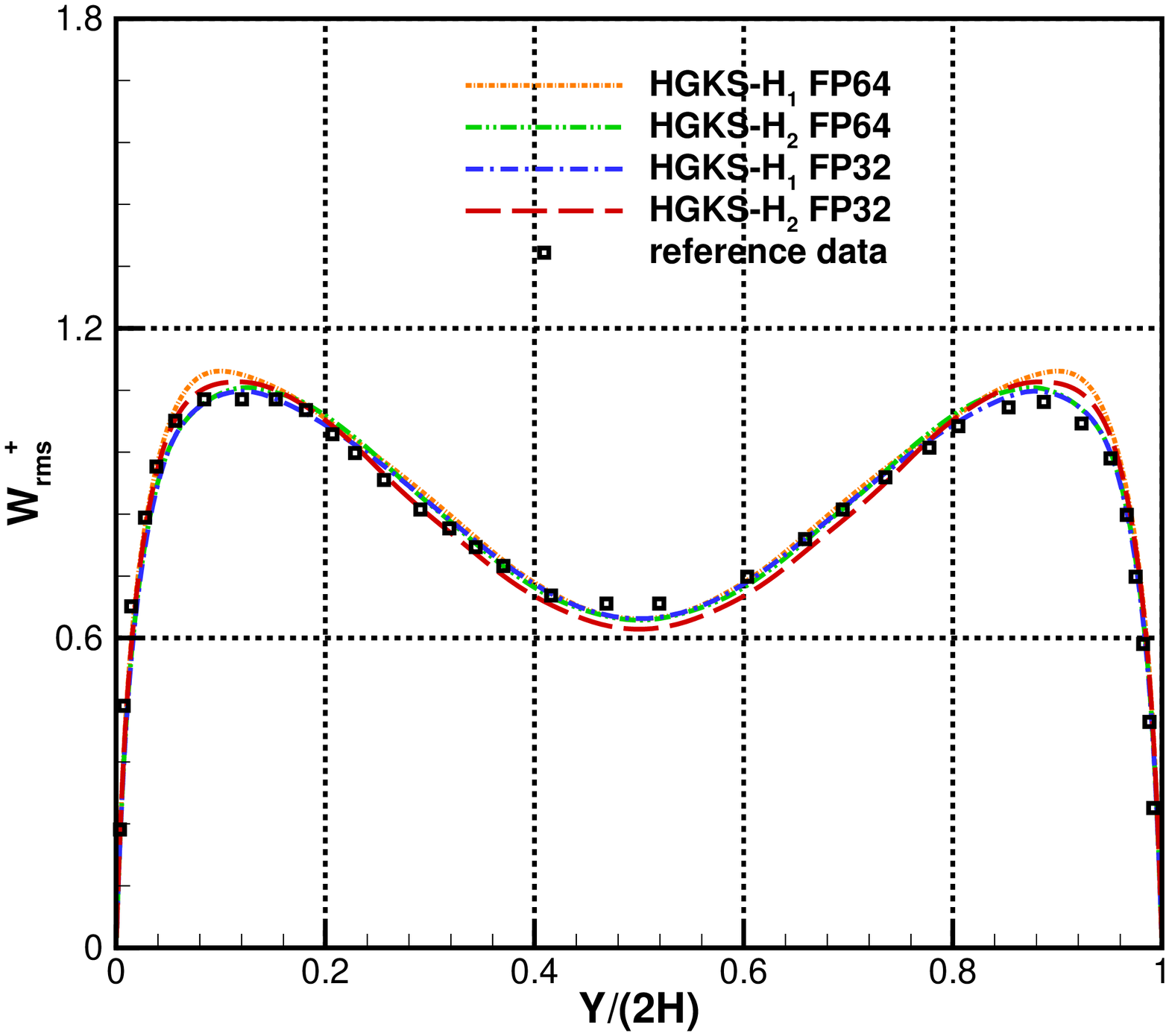}
\includegraphics[width=0.45\textwidth]{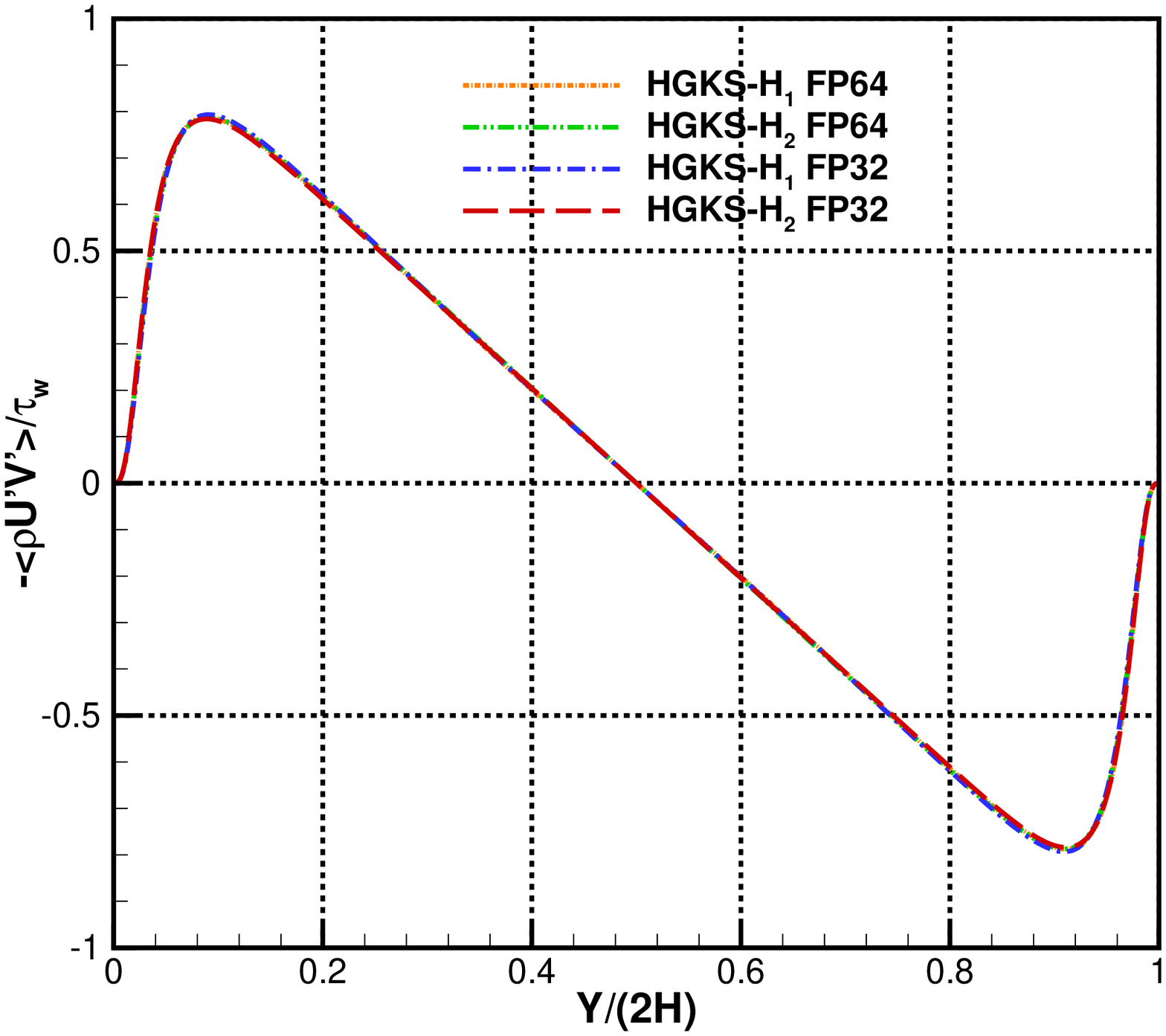}
\includegraphics[width=0.45\textwidth]{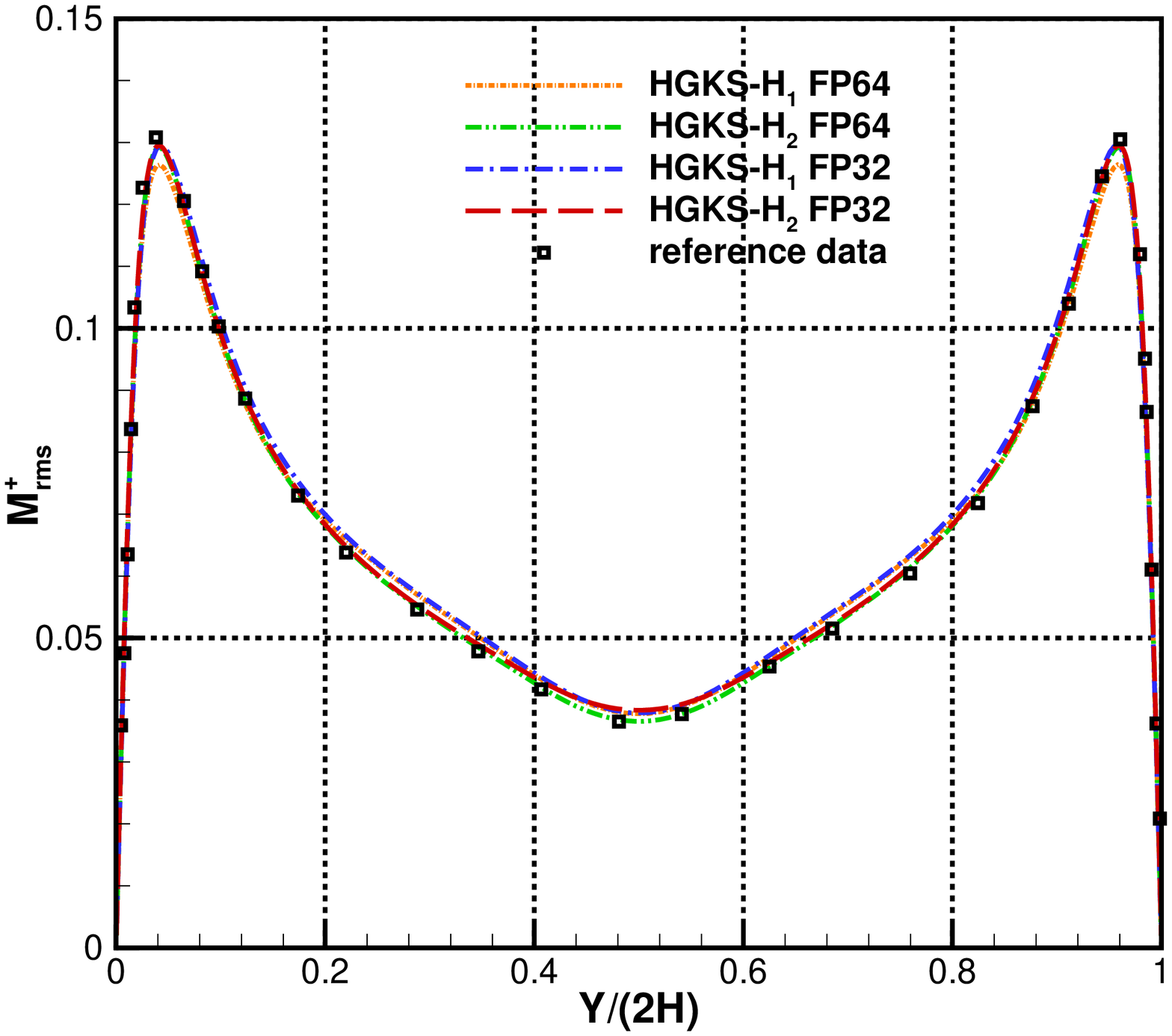}
\includegraphics[width=0.45\textwidth]{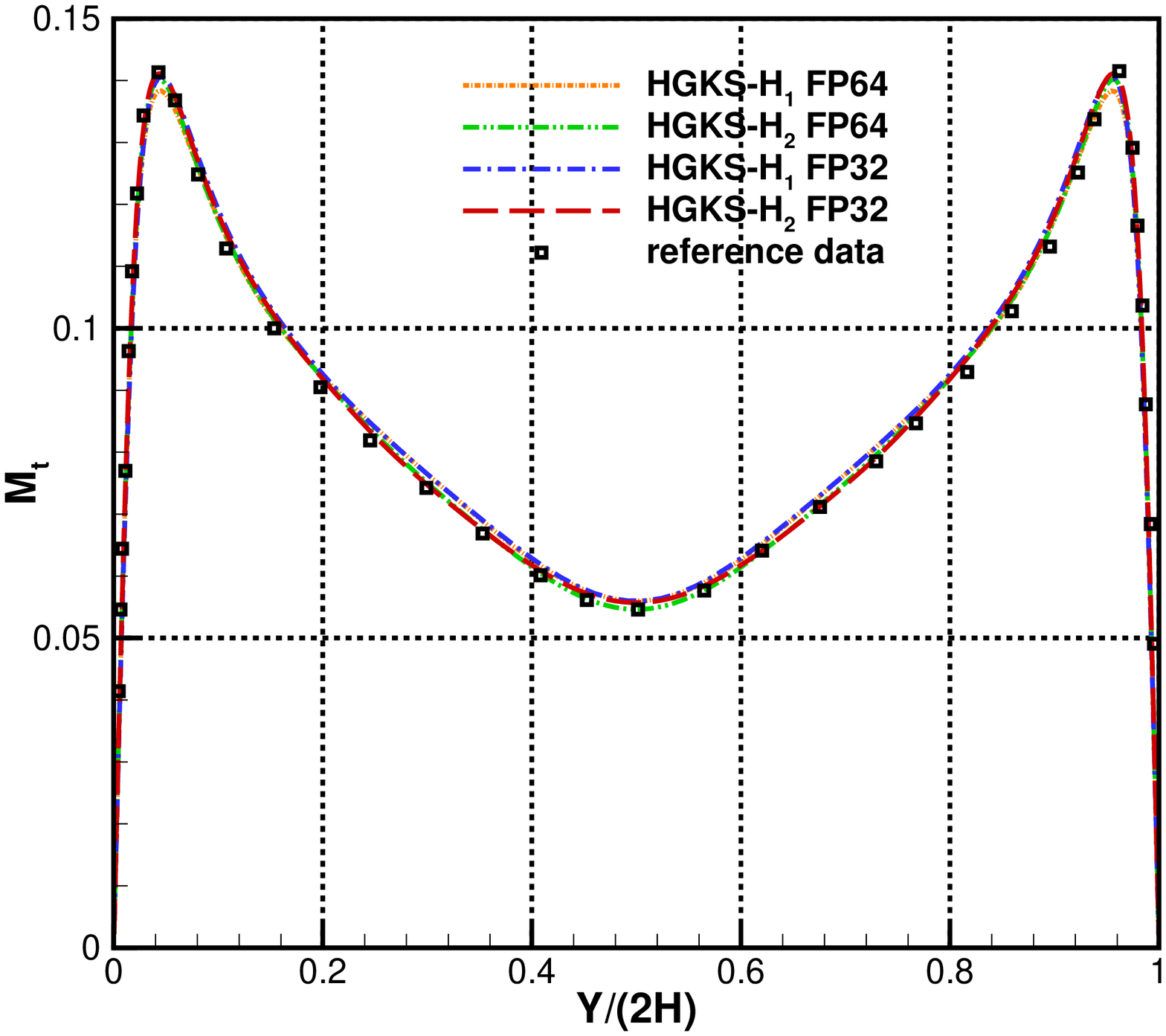}
\caption{\label{channel_fluctuation_2} Turbulent channel flow: the
root-mean-square fluctuation velocity $U_{rms}^+$, $V_{rms}^+$,
$W_{rms}^+$, Reynolds stress $\langle-\rho U^{'} V^{'}
\rangle/\langle \tau_{w} \rangle$, the root-mean-square of Mach
number $M_{rms}^+$ and turbulent Mach number $M_t$ profiles for
compressible turbulent flow. The reference data is given in
\cite{DNS-Li}.}
\end{figure}

To validate the accuracy of HGKS and address the effect with
different precision, the statistical turbulent quantities are
provided for the nearly incompressible and compressible turbulent
channel flows. For FP32-FP64 comparison, all cases restarted from
the same flow fields, and $200 H/U_c$ statistical periods are used
to average the computational results for all cases. For the nearly
incompressible channel flow, the logarithmic formulation is given by
\begin{align}\label{log-law}
U^+ = \frac{1}{\kappa}\ln Y^+ + B,
\end{align}
where the von Karman constant $\kappa=0.40$ and $B=5.5$ is given for
the low Reynolds number turbulent channel flow
\cite{kim1987turbulence}. The mean streamwise velocity profiles with
a log-linear plot are given in Fig.\ref{channel_ave_1}, where the
HGKS result is in reasonable agreement with the spectral results
\cite{kim1999turbulence}. In order to account for the mean property
of variations caused by compressibility, the Van Driest (VD)
transformation \cite{vdtransformation} for the density-weighted
velocity is considered
\begin{align}\label{log-law-vd}
{\left\langle U \right\rangle}_{VD}^+ = \int_{0}^{{\left\langle U \right\rangle}^+} \large \big(\frac{\left\langle \rho \right\rangle}{\left\langle \rho_w \right\rangle} \big)^{1/2} \text{d} {\left\langle U \right\rangle}^+,
\end{align}
where $\langle\cdot\rangle$ represents the mean average over
statistical periods and the X- and Z-directions. For the
compressible flow, the transformed velocity is expected to satisfy
the incompressible log law Eq.\eqref{log-law}. The mean streamwise
velocity profiles $\langle U \rangle^+$ and $\langle U
\rangle_{VD}^+$ with VD transformation as Eq.\eqref{log-law-vd} are
given in Fig.\ref{channel_ave_2} in log-linear plot. HGKS result is
in reasonable agreement with the reference DNS solutions with the
mesh refinements \cite{DNS-Li}. The mean velocity with FP32 and FP64
precision are also given in Fig.\ref{channel_ave_1} and
Fig.\ref{channel_ave_2}. It can be observed that the qualitative
differences between different precision are negligible for
statistical turbulent quantities with long time averaging.

For the nearly incompressible turbulent channel flow, the time
averaged normalized root-mean-square fluctuation velocity profile
$U_{rms}^+$, $V_{rms}^+$, $W_{rms}^+$ and normalized Reynolds stress
profiles $-<U'V'>$ are given in Fig.\ref{channel_fluctuation_1} for
the cases $G_1$ and $G_2$. The root mean square is defined as
$\phi_{rms}^{+} = \sqrt{(\phi -\langle \phi \rangle)^2}$ and
$\phi'=\phi -\langle \phi \rangle$, where $\phi$ represents the flow
variables. With the refinement of mesh, the HGKS results are in
reasonable agreement with the reference date, which given by the
spectral method with $256 \times 193 \times 192$ cells
\cite{kim1999turbulence}. It should also be noted that spectral
method is for the exact incompressible flow, which provides accurate
incompressible solution. The HGKS is more general for both nearly
incompressible flows and compressible flows. For the compressible
flow, the time averaged  normalized root-mean-square fluctuation
velocity profiles $U_{rms}^+$, $V_{rms}^+$, $W_{rms}^+$, Reynolds
stress $\langle-\rho U^{'} V^{'} \rangle/\langle \tau_{w} \rangle$,
the root-mean-square of Mach number $M_{rms}^+$ and turbulent Mach
number $M_t$  are presented in Fig.\ref{channel_fluctuation_2} for
the cases $H_1$ and $H_2$. The turbulent Mach number is defined as
$M_t = q/\left\langle c \right \rangle$, where $q^2 =\langle
(U^{'})^2+(V^{'})^2+(W^{'})^2\rangle$  and $c$ is the local sound
speed. The results with $H_2$ agree well with the refereed DNS
solutions, confirming the high-accuracy of HGKS for DNS in
compressible wall-bounded turbulent flows. The statistical
fluctuation quantities profiles with FP32 precision and FP64
precision are also given in Fig.\ref{channel_fluctuation_1} and
Fig.\ref{channel_fluctuation_2}. Despite the deviation of
$U_{rms}^+$ for different precision, the qualitative differences in
accuracy between FP32 and FP64 precision are acceptable for
statistical turbulent quantities with such a long time average.

\begin{figure}[!h]
\centering
\includegraphics[width=0.45\textwidth]{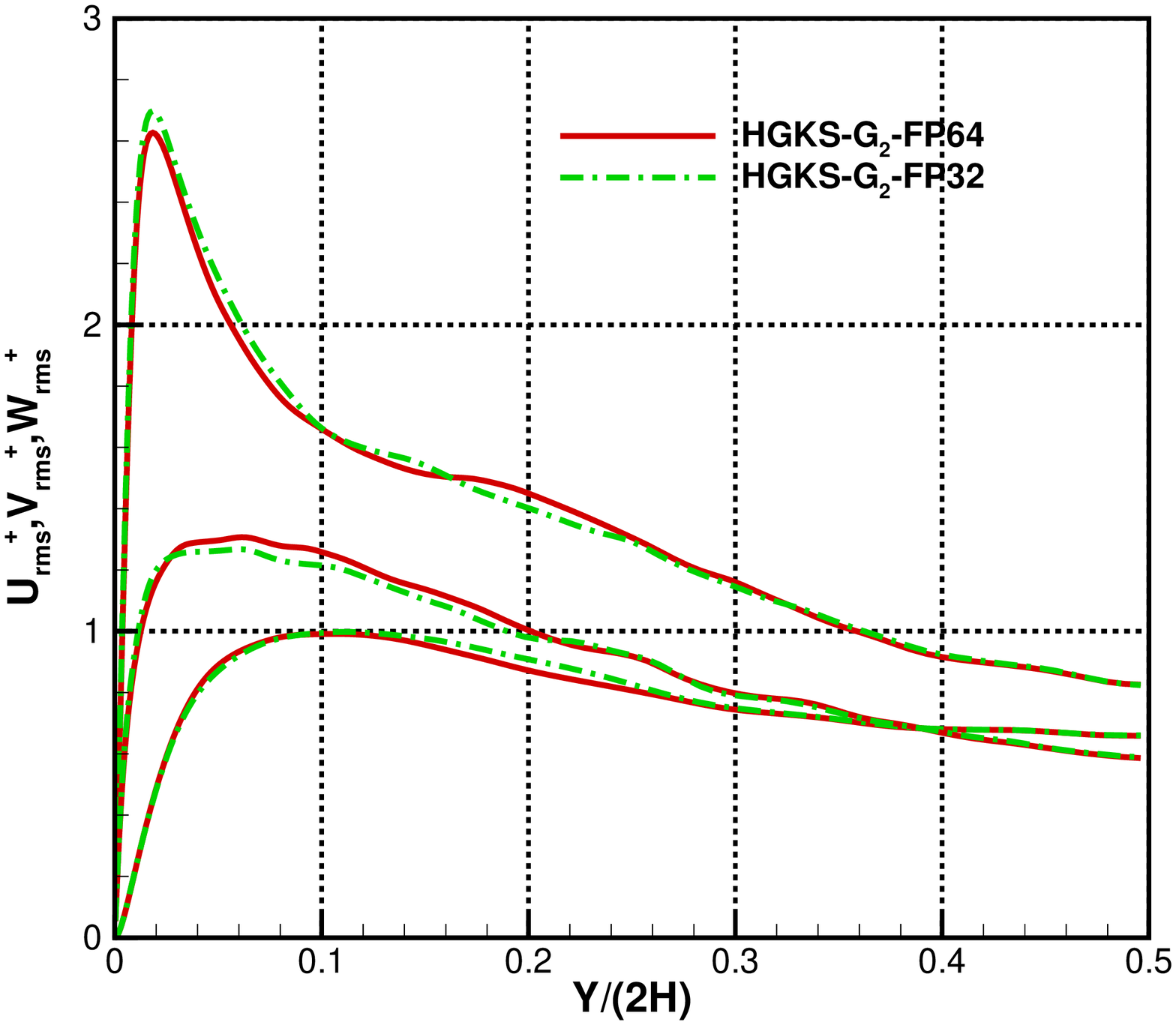}
\includegraphics[width=0.45\textwidth]{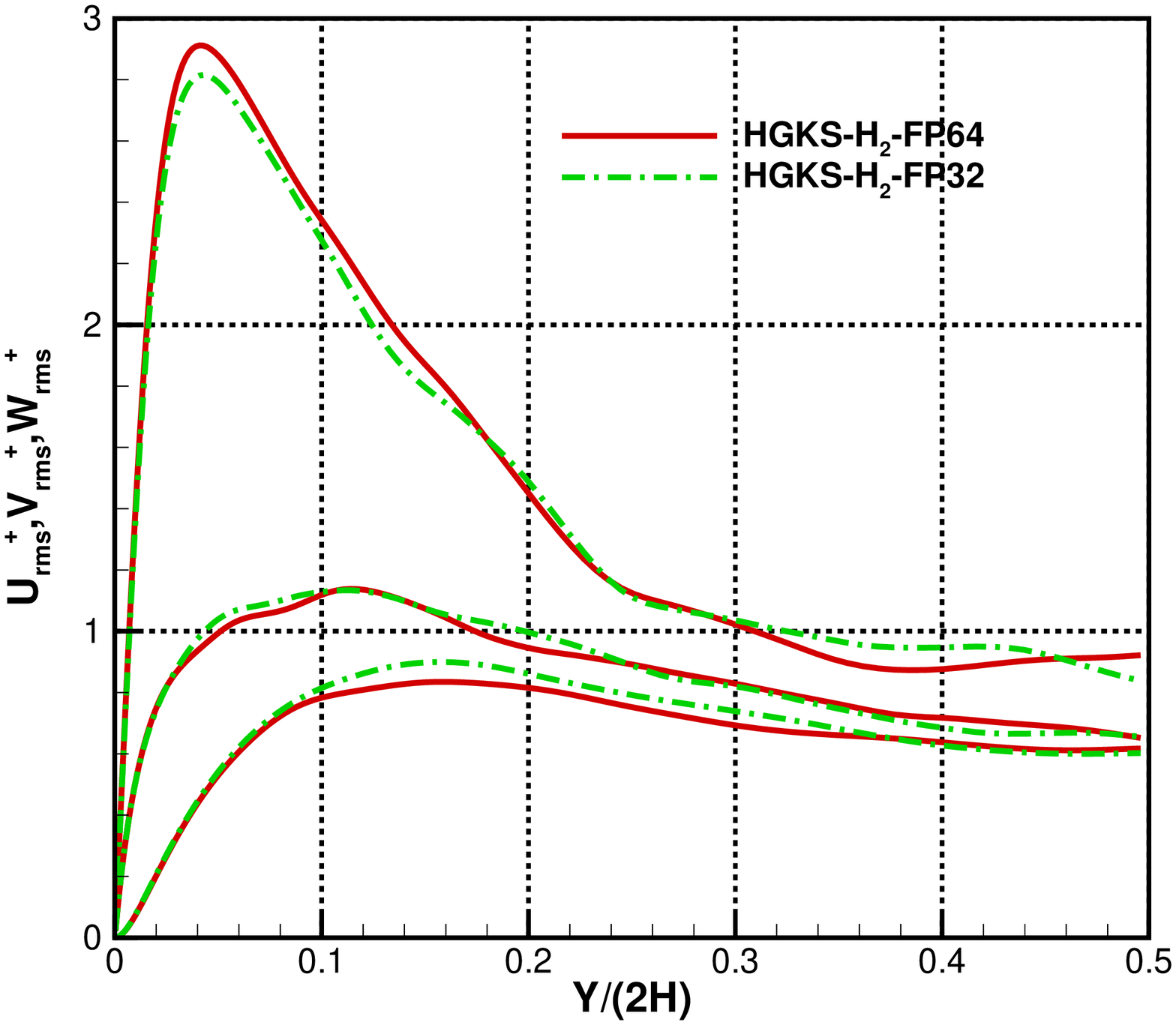}
\includegraphics[width=0.45\textwidth]{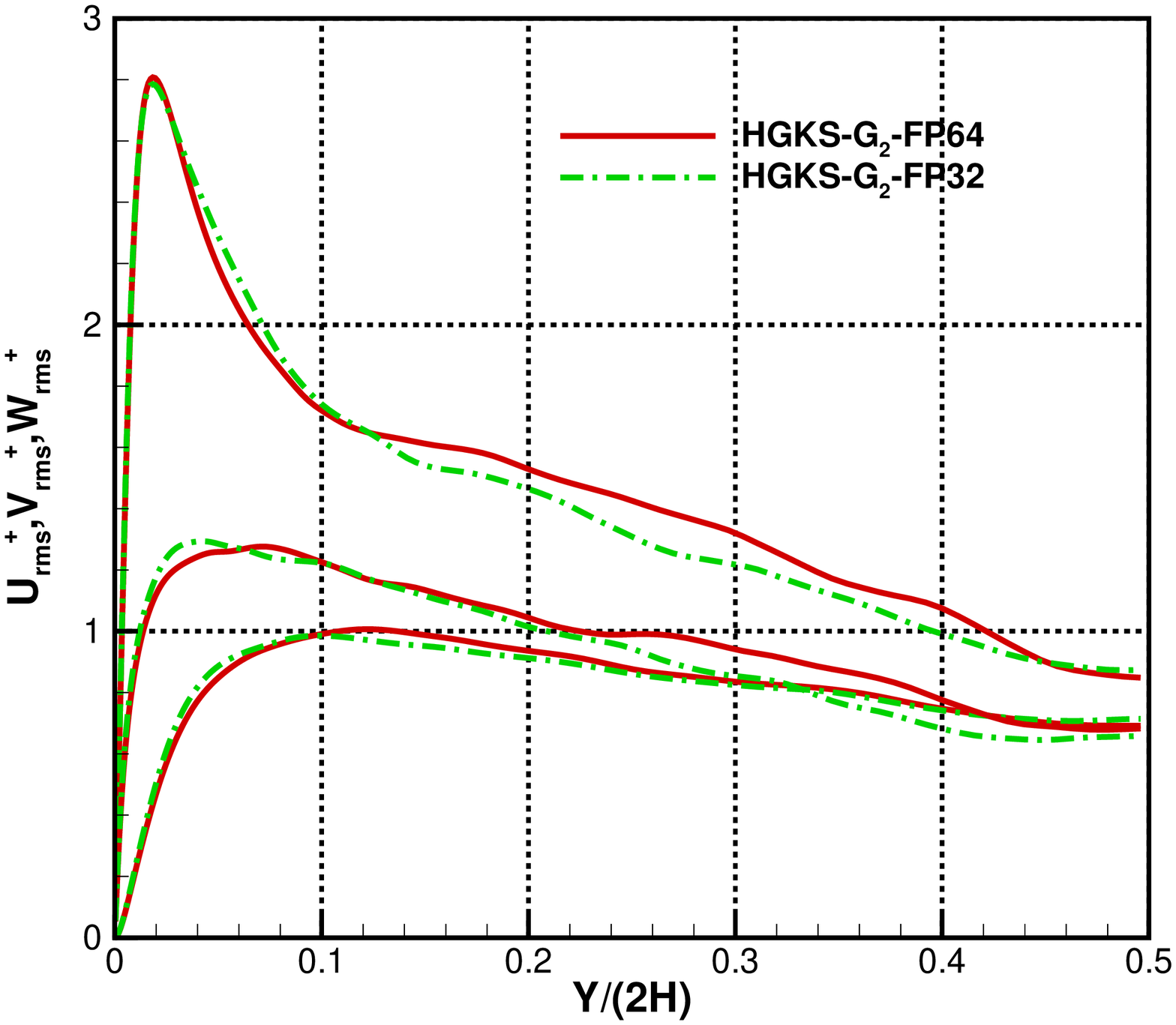}
\includegraphics[width=0.45\textwidth]{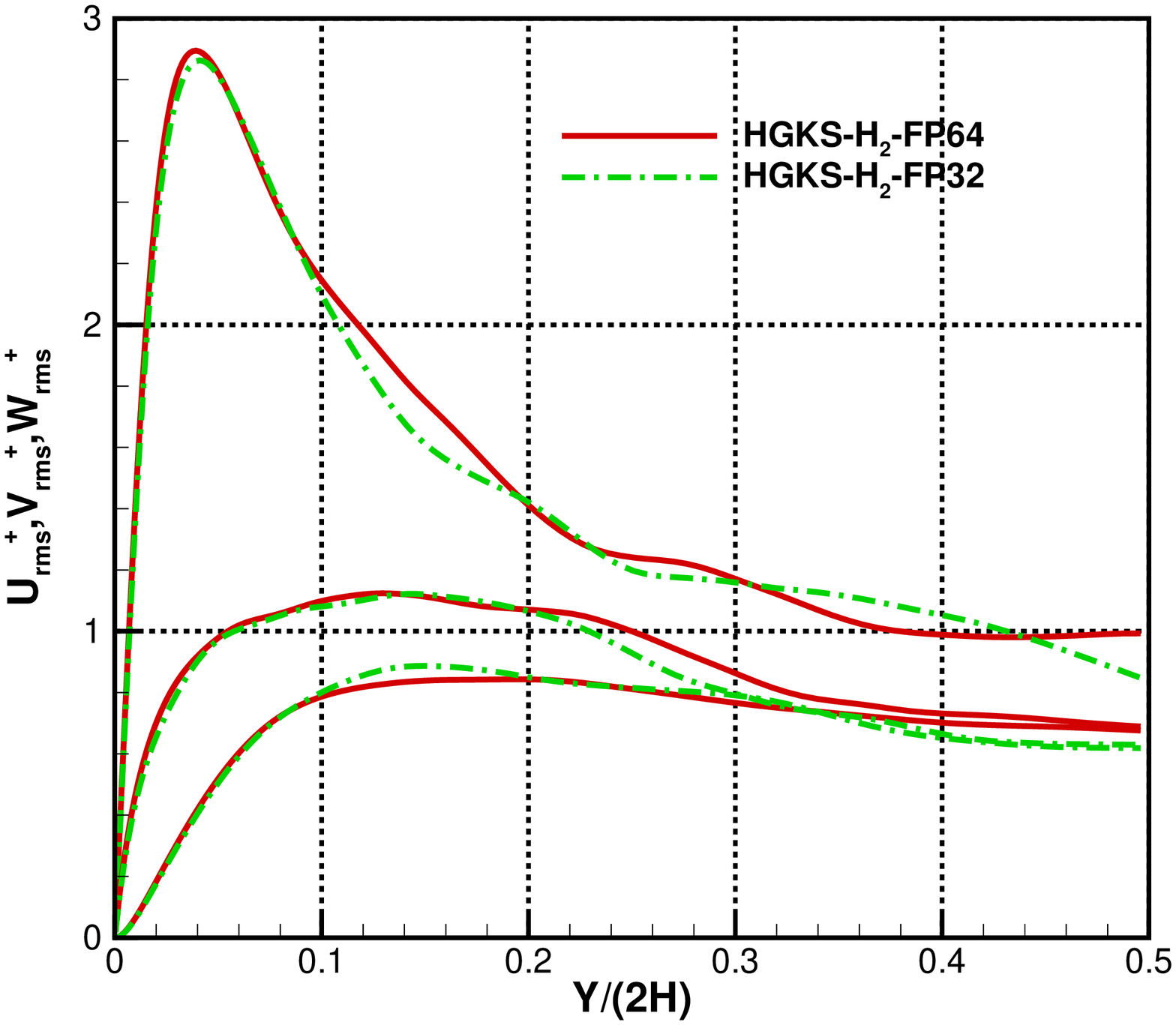}
\caption{\label{channel_fluctuation_3} Turbulent channel flow: the
instantaneous profiles of $U^+$, $V^+$, $W^+$ for $G_2$ and $H_2$ at
$t=20H/U_c$ (top) and $200H/U_c$ (bottom) with FP32 and FP64
precision.}
\end{figure}

The instantaneous profiles of $U^+_{rms}$, $V^+_{rms}$, $W^+_{rms}$
at $t=20H/U_c$ and $200H/U_c$ are shown in
Fig.\ref{channel_fluctuation_3} for the cases $G_2$ and $H_2$ with
FP32 and FP64 precision, where cases restarted with the same initial
flow field for each case. Compared with the time averaged profiles,
the obvious deviations are observed for the instantaneous profiles.
The deviations are caused by the round-off error of FP32 precision,
which is approximately equals to or larger than the  errors of
numerical scheme. Since detailed effect of the round-off error is
not easy to be analyzed, FP32 precision may not be safe for DNS in
turbulence especially for time-evolutionary turbulent flows. In
terms of the problems without very strict requirements in accuracy,
such as large eddy simulation (LES) and Reynolds-averaged
Navier–Stokes (RANS) simulation in turbulence, FP32 precision may
be used due to its improvement of efficiency and reduction of
memory. It also strongly suggests that the FP64 precision
performance of GPU still requires to be improved to accommodate the
increasing requirements of GPU-based HPC.

\section{Conclusion}
Based on the multi-scale physical transport and the coupled
temporal-spatial gas evolution, the HGKS provides a workable tool
for the numerical study of compressible turbulent flows. In this
paper, to efficiently conduct large-scale DNS of turbulence, the
HGKS code is developed with single GPU using CUDA architecture, and
multiple GPUs using MPI and CUDA. The multiple GPUs are distributed
across multiple CPUs at the cost of having to coordinate network
communication via MPI. The Taylor-Green vortex problems and
turbulent channel flows are presented to validate the performance of
HGKS with multiple Nvidia TITAN RTX and Nvidia Tesla V100 GPUs. We
mainly concentrate on the computational efficiency with single GPU
and multiple GPUs, and the comparisons between FP32 precision and
FP64 precision of GPU. For single-GPU computation, compared with the
OpenMP CPU code using Intel Core i7-9700, 7x speedup is achieved by
TITAN RTX and 16x speedup is achieved with Tesla V100. The
computational time of single Tesla V100 GPU is comparable with the
MPI code using $300$ supercomputer cores with Intel Xeon E5-2692.
For multiple GPUs, the HGKS code scales properly with the number of
GPU. It can be inferred that the efficiency of GPU code with $8$
Tesla V100 GPUs approximately equals to that of MPI code with $3000$
CPU cores. Compared with FP64 precision simulation, the efficiency of HGKS can be
improved and the memory is reduced with FP32 precision. However,
with the long time computation in compressible turbulent channels
flows, the differences in accuracy appear. Especially, the
instantaneous statistical turbulent quantities is not acceptable
using FP32 precision. The choice of precision should depend on the
requirement of accuracy and the available computational resources.
In the future, more challenging compressible flow problems, such as
the supersonic turbulent boundary layer and the interaction of shock
and turbulent boundary layer \cite{adams2000direct,wu2007direct},
will be investigated with efficient multiple-GPU accelerated HGKS.

\section*{Ackonwledgement}
This research is supported by National Natural Science Foundation of
China (11701038), the Fundamental Research Funds for the Central
Universities. We would thank Prof. Shucheng Pan of Northwestern
Polytechnical University, Dr. Jun Peng and Dr. Shan Tang of Biren
Technology for insightful discussions in GPU-based HPC.

\end{document}